\documentclass[makeidx]{crmp-l}

\makeindex

\usepackage{amsmath}
\usepackage{amssymb}
\usepackage{epsfig}

\newtheorem{theorem}{Theorem}[section]
\newtheorem{lemma}[theorem]{Lemma}
\newtheorem{prop}[theorem]{Proposition}

\theoremstyle{definition}

\theoremstyle{remark}

\numberwithin{equation}{section}

 \newcommand{\bs}[1]{\boldsymbol{#1}}

 \newcommand{\NN}{\mathbb N}
 \newcommand{\RR}{\mathbb R}
 \newcommand{\CC}{\mathbb C}
 \newcommand{\ZZ}{\mathbb Z}
 \newcommand{\QQ}{\mathbb Q}

 \newcommand{\cA}{\mathcal{A}}
 \newcommand{\cL}{\mathcal{L}}
 
 \newcommand{\cK}{\mathcal{K}}
 \newcommand{\cM}{\mathcal{M}}
 \newcommand{\cP}{\mathcal{P}}
 \newcommand{\cF}{\mathcal{F}}
 \newcommand{\cG}{\mathcal{G}}

 \newcommand{\conv}{\mbox{\Huge \raisebox{-0.2ex}{$\ast$}}}
 \newcommand{\btimes}{\mbox{\huge \raisebox{-0.2ex}{$\times$}}}

\begin{document}

\title{Self-Similar Measures for Quasicrystals}

\author{Michael Baake}

\address{Institut f\"ur Theoretische Physik, Universit\"at T\"ubingen,
         Auf der Morgenstelle 14, 72076 T\"ubingen, Germany}

\email{michael.baake@uni-tuebingen.de}
\email{baake@miles.math.ualberta.ca}

\author{Robert V.\ Moody}
\address{Department of Mathematical Sciences, University of Alberta,
         Edmonton, Alberta, Canada T6G 2G1}
\email{rvm@miles.math.ualberta.ca}

%    General info
% \subjclass{Primary 54C40, 14E20; Secondary 46E25, 20C20}
% \date{January 1, 1994 and, in revised form, June 22, 1994.}

\dedicatory{Dedicated to Peter A.\thinspace B.\thinspace Pleasants on
the occasion of his 60th birthday.}

\thanks{This work was supported by the Natural Sciences and Engineering
Research Council of Canada (NSERC) and by the German Research Council (DFG)}

\begin{abstract}
We study self-similar measures of Hutchinson type, defined by compact families 
of contractions, both in a single and multi-component setting. The results are 
applied in the context of general model sets to infer, via a generalized 
version of Weyl's Theorem on uniform distribution, the existence of invariant 
measures for families of self-similarities of regular model sets.
\end{abstract}

\maketitle

\section*{Introduction}

There have been two very successful approaches to generate aperiodic
sets with features of long-range internal order. The first is by
creating tilings by the method of inflation followed by decomposition
using a finite set of proto-tiles. The second is by creating point sets
through the method of cut and project sets, or model sets as we call
them here. Neither theory subsumes the other and they both have their
own particular virtues. However, they have a considerable overlap. It is
easy to replace a tiling by an equivalent set of points (by selecting
suitable points from each type of tile) and in many cases the result is
a model set or, more generally, a union of several model sets, one for 
each type of point.  Conversely, there are many ways to obtain a tiling 
from a point set (for instance by using the Voronoi cells, or the dual 
Delone cells determined by them), and the equivalence concept of mutual
local derivability, see \cite{MB,P} and references given there, is an
adequate tool to make this connection precise. So it
is natural to study the most notable feature of inflation tilings,
namely their self-similarity, in the context of model sets, and indeed,
even when no simple tiling is in sight, many interesting model sets have
striking self-similarity.

The objective of this paper is to set up some of the machinery that
makes such a study possible and to show how naturally it can be
associated with families of self-similar measures on locally compact
Abelian groups.

There is one initial hurdle which is not usually considered in the
study of model sets. In order to have any sort of reasonable
correspondence with the tiling world, and to create a useful theory, we
need to have a {\em multi-component} model set 
\index{model set!multi-component} context in which there are a finite
number of model sets, all based on the same cut and project scheme,
that are mutually coupled by the self-similarities. After all, almost
all inflation tiling systems have several types of tiles, and the
decomposition of inflated tiles typically involves all of these
various types simultaneously.

Thus the situation that we envision consists of a {\em family} of model 
sets \index{family of model sets}
$\Lambda_1, \dots ,\Lambda_n$ in some real space $\RR^m$, all based on
the same cut and project scheme, and a set of families of inflationary
mappings $F_{ij}$ on $\RR^m$ with the property that 
\begin{equation}\label{inflationEquationIntro} 
    \Lambda_i \; = \; \bigcup_{j=1}^n \bigcup_{f \in F_{ij}} f(\Lambda_j)\, , 
    \qquad 1 \le i \le n \,.  
\end{equation}

There are two new features that are different from the tiling
situation. In inflation tilings, the usual idea is that the tiles
only overlap on their boundaries which are of measure $0$. That would
be equivalent to some form of disjointness in Eq.\
(\ref{inflationEquationIntro}) which we do {\em not} wish to assume. Secondly,
in tilings, the sets $F_{ij}$ of mappings are assumed finite. Again, we
do {\em not} wish to make this assumption. In fact, it is quite useful to
make $F_{ij}$ consist of {\em all} possible mappings of a certain type (for
example all affine mappings) that are consistent with Eq.\
(\ref{inflationEquationIntro}), and this is in general an infinite set.

As soon as overlapping is allowed, there naturally arises the question
of whether there is implicit in Eq.~(\ref{inflationEquationIntro}) a
corresponding set of relative weights which measure the ``frequency"
of occurrence of the points of $\Lambda_i$ coming from the substitution
process (\ref{inflationEquationIntro}). If we were counting, this
could be loosely construed as noting the occurrences of points of
$\Lambda_i$ with their multiplicities as they appear in the righthand
side of (\ref{inflationEquationIntro}). This links up to a recent
approach of Lagarias and  Wang \cite{LW}. However, in our situation, there is no
reason to assume that our weights can be normalized to be integral
nor that there are only a finite number of contributions to
$\Lambda_i$ in the right-hand side of (\ref{inflationEquationIntro}).

This then is the primary goal of this contribution: to discuss the existence
and nature of self-similar densities on systems of model sets which are
coupled together by a substitution system of type
(\ref{inflationEquationIntro}). Our definition of model
sets is based on arbitrary locally compact Abelian groups as internal
spaces, and so is more flexible and handles more situations (see for
instance \cite{BMS,LM}) than the usual method of projection from Euclidean
spaces. This extra generality requires a little more care than usual
but it is remarkable how much of the theory is natural in this context.
 
The method by which we attack the
problem is to use the formalism of the underlying cut and project
scheme to pass everything over to the internal side of the picture,
i.e.\ onto the locally compact Abelian group that is controlling the
projection. The advantage of doing this is that the system of
inflationary mappings turns into a family of contractions, and what has
started off as a problem in the domain of discrete mathematics turns
into one of analysis. A primary virtue of systems of contractions is
the ready-built Hutchinson theory of iterated function systems
with their attractors and self-similar measures. 
In our multi-component setting with infinite families of
mappings, we need a slight variation on this theme, which
occupies the first five Sections of the paper. These parts of the paper
have nothing in particular to do with model sets, but rather are a
development of Hutchinson's theory in the multi-component situation
where the coupling is by {\em compact} families of contractions 
\index{compact sets of contractions} each of which
has its own, essentially arbitrarily pre-given positive Borel measure.
An important part of this is determining some conditions under which
the self-similar measures are in fact absolutely continuous and, more
importantly, when the representing $L^1$-functions, the Radon-Nikodym densities,
are actually continuous (continuously representable self-similar measures), 
for it is only then that we can bring back information to the discrete side again.

Section \ref{model-weyl} of this article brings in the model sets and develops
the  mathematics that allows us to pass information back from the
internal side to the model set side in $\RR^m$.  The primary tool here
is Weyl's theory of uniform projection, but we have to redevelop
this in the context of locally compact Abelian groups and model sets.

In Section \ref{model-dens} we are, at last, set to tackle the problem of
determining self-similar distribution of weights (also called self-similar
densities from now on) on the model sets themselves. The
Weyl theory applies only to continuous functions, so we can only refer 
to it when our self-similar measures have continuously representable
Radon-Nikodym densities.
Fortunately, this is the case in a number of interesting situations. We
provide a general description of the situations in which such
self-similar densities exist, and then, in Section \ref{sec-examples}, we offer a
number of examples which illustrate what we have achieved.

The results obtained here extend those of two previous
papers \cite{BM,BMmulti}. There, we considered only model
sets based on Euclidean internal spaces, and the context was not primarily 
measure-theoretical as it is here. The method of dealing
with multi-component model sets also differed from the product
approach that we have adopted in this article. There are nonetheless several
examples of self-similar densities (called invariant densities
there) in those papers which the reader may find of interest.

\section{Compact families of contractions and attractors}
\label{comp-fam}

Let us first review some basic facts from the theory of iterated function
systems, both finite and compact, in a way that is adequate for our needs.

\subsection{Hutchinson's contraction principle}

Let $X$ be a complete metric space with metric $d$. 
We denote by $\cK X$ the set of all non-empty compact subsets of $X$. Let
$d(x,U):=\inf\{d(x,u)\mid u\in U\}$ be the distance of $x$ from $U$.
Note that $d(x,U)=0$ implies $x\in\overline{U}$.
For $\varepsilon >0$, the $\varepsilon$-$fringe$\index{fringe}
of a subset $U \in \cK X $ is 
\begin{equation}
    [U]_\varepsilon \; := \; \{ x \in X \mid d(x,U) < \varepsilon \} \, .
\end{equation}
In view of the set $\cK X$, we introduce the {\em Hausdorff metric}: 
\index{Hausdorff metric} for $U, V \in \cK X$, it is defined by
\begin{equation} 
   d_H(U,V) \; := \;
     \inf\{ \varepsilon >0 \mid U \subset [V]_\varepsilon 
               \mbox{ and } V \subset [U]_\varepsilon  \} \, .
\end{equation}
Note that for singletons $U = \{u\}$ and $V=\{v\}$, one has
$d_H(\{u\},\{v\}) = d(u,v)$. An alternative way to determine $d_H(U,V)$ is
$$ d_H(U,V) \; = \;
   \sup\{ d(u,V), \, d(v,U)\mid u\in U, \, v\in V \}\, .$$

Relative to the Hausdorff metric, $\cK X$ is again a complete metric space.
If now $U_j,V_j \in \cK X, \; j \in J$, are two sets of compact subsets of $X$ 
then it is easy to see \cite[Note 2.1.6]{Wicks} that 
\begin{equation} \label{unionConst}
    d_H \Bigl(\,\mbox{$\bigcup_{j\in J}$} \, U_j, 
                \mbox{$\bigcup_{j\in J}$} \, V_j \Bigr) 
             \; \le \; \sup_{j\in J} \, d_H(U_j,V_j) \, .
\end{equation}

Given complete metric spaces $X,Y$, we consider the space $C(X,Y)$
of all continuous mappings of $X$ into $Y$, endowed with the compact-open
topology. This topology has the property of making the evaluation 
\index{evaluation} maps
\begin{equation} \label{evaluation}
    {\rm eval}_x: \; C(X,Y) \longrightarrow Y \; ; \quad f \mapsto f(x)
\end{equation}
continuous \cite[Thm.\ 7.4]{Kelley}. 
There is a natural extension of mappings that leads to
\begin{equation} \label{extension}
    \cK(\centerdot) : \; C(X,Y) \longrightarrow C(\cK X, \cK Y)
\end{equation}
with $(\cK (f)) (U) := f(U)$.
This mapping is continuous \cite[Prop.\ 2.5.1]{Wicks}\footnote{Wicks' book 
\cite{Wicks} is a 
great place to look for information on spaces of compact sets, and it was an 
important source for the part of this paper on compact spaces of maps. 
Unfortunately, or fortunately, according to one's taste, Wicks uses 
non-standard analysis to streamline his presentation, so `standard' readers 
need to look elsewhere for the proofs, or to adapt them.}.

A mapping $f:X \longrightarrow Y$ of metric spaces is {\em Lipschitz} if 
there is \index{Lipschitz function} an $r>0$ with 
$d_Y(f(x^{}_1),f(x^{}_2)) \le r \, d_X(x^{}_1,x^{}_2)$ for all 
$x^{}_1,x^{}_2 \in X$. If $f$ is Lipschitz then the infimum of all such $r$ 
is the {\em Lipschitz constant} $r^{}_f$ \index{Lipschitz constant} of $f$. 
If $r^{}_f <1$ then $f$ is called a {\em contraction}. \index{contraction}
The set of all Lipschitz maps from $X$ to $Y$ with Lipschitz constant 
equal to $r$ (resp.\ at most $r$) is denoted by ${\rm Lip}(r,X,Y)$ 
(resp.\ ${\rm Lip}(\le\! r,X,Y)$). Lipschitz functions are clearly uniformly 
continuous mappings. We observe:
\begin{equation} \label{lipK}
    f \in {\rm Lip}(r,X,Y) \quad \Longrightarrow \quad 
    \cK (f) \in {\rm Lip}(r, \cK X, \cK Y) \, . 
\end{equation}
Note in particular that $r^{}_f = r^{}_{\cK (f)}$, as is seen
by looking at singleton sets.

We are interested in {\em union maps}. \index{union maps}
Let $X,Y$ be complete metric spaces and 
let $F = \{f_1, f_2, \dots , f_N\}$ be a set of continuous maps from $X$ to 
$Y$. We define two mappings, both given the same name:
\begin{equation} \label{unionMap}
\begin{array}{crlcl}
   {}_{}^\cup\! F: & X \rightarrow \cK Y \, ; &
         {}_{}^\cup\! F(x) & := & \mbox{$\bigcup_{i=1}^N$} \{f_i(x)\} \\
  \vphantom{\mbox{\Huge F}}
   {}_{}^\cup\! F: & \cK X \rightarrow \cK Y \, ; & 
         {}_{}^\cup\! F(U) & := & \mbox{$\bigcup_{i=1}^N$} \; f_i(U) \, .
\end{array}
\end{equation}
In view of (\ref{unionConst}), the union map ${}_{}^\cup\! F$ is Lipschitz
if all mappings in $F$ are Lipschitz, and we have
$r^{}_{{}_{}^\cup\! F} \le \sup \{r^{}_{f} \mid f\in F\}$.

When $X = Y$ and $F$ consists of contractions, then $F$ is called an 
{\em iterated function system} (IFS). \index{iterated function system}
The principal result, see \cite[{\S}  3.1]{Hutchinson} or
\cite[Prop.\ 3.1.1]{Wicks}, is based on the general
Banach contraction principle and reads as follows.
\begin{theorem} {\rm (Hutchinson's contraction principle)} 
   Let $F$ be an {\rm IFS} on a complete metric space $X$. 
   Then there is a unique 
   $W \in \cK X$ which is a fixed point of ${}_{}^\cup\! F$, i.e.\  
   $W = \mbox{$\bigcup_{f\in F}$} \, f(W)$. Furthermore, for any 
   $Z\in \cK X$, $({}_{}^\cup\! F)^{\ell}(Z)$ converges to $W\in\cK X$
   in the Hausdorff metric, as $\ell\to\infty$.   \qed
\end{theorem}
$W$ is the {\em attractor} \index{attractor} of the IFS.  
Our aim, in this and in the following Section, is to generalize this result 
in two directions:
\begin{enumerate}
  \item to compact sets of contractions (which is also well known);
  \item to products of metric spaces (in what we call the multi-component 
               situation).
\end{enumerate}

\subsection{Compact sets of contractions}

Let $X,Y$ be complete metric spaces. We form $\cK X$ as above.
We remark that the union of any compact subset of $\cK X$ is a
compact subset of $X$, i.e.\
\begin{equation} \label{compactUnion}
     C\in \cK \cK X \quad \Longrightarrow \quad 
      \mbox{$\bigcup$}C \, := \,
      \mbox{$\bigcup_{U \in C}$} U  \;\in\; \cK X \, ,
\end{equation}
and that the union map $\bigcup\! : \cK\cK X \to \cK X$ is continuous, 
cf.~\cite[Ch.\ 1.5]{Wicks}.

Consider the space $C(X, \cK Y)$ equipped with the compact-open topology. 
Let $F \in \cK C(X,\cK Y)$, i.e.\ $F$ is a compact subset of contiuous 
mappings from $X$ into the space of compact subsets of $Y$. In view of 
(\ref{evaluation}), we have, for all $x \in X$,
\begin{equation}
    F(x) \; := \; \{f(x) \mid f \in F \}  \; \in \;  \cK \cK Y \, , 
\end{equation}
i.e.\ $F(x)$ is compact. From (\ref{compactUnion}), we deduce that 
$\bigcup_{f\in F} f(x)$ is also compact. Thus, as in (\ref{unionMap}), 
we have a new mapping
\begin{equation} \label{compactUnionMap}
\begin{array}{rl}
{}_{}^\cup\! F : & X   \longrightarrow \cK Y \\
                 & \,x \;\,\mapsto \,\;\mbox{$\bigcup_{f\in F}$}\, f(x) \; .
\end{array}
\end{equation}

With these preliminary definitions out of the way, we define $G$ to be a 
{\em compact admissible family of Lipschitz mappings} \cite[Sec.\ 3.1]{Wicks} 
\index{admissible Lipschitz mappings} from  $X$ to $Y$ if 
\begin{itemize}
\item[{\bf C1}]  $G $ is a compact set of Lipschitz mappings from $X$ to $Y$;
\item[{\bf C2}]  there exists $r > 0$ such that $r^{}_g\le r$ for all $g\in G$.
\end{itemize}
If $r <1$ then we call $G$ a {\em compact admissible family of contractions}.

Let $G$ be a compact admissible family of Lipschitz  mappings from $X$ to $Y$. 
{}From Eq.\ (\ref{extension}), $\cK (G) \subset C(\cK X, \cK Y)$ is compact, 
and since the $\cK$ operator preserves Lipschitz constants, $\cK (G)$ is itself 
a compact admissible family of Lipschitz mappings.  Thus, from 
(\ref{compactUnionMap}) we obtain
\begin{equation} \label{cKMap}
   {}_{}^\cup\! \cK G := {}_{}^\cup\! (\cK (G))\; \in \; C(\cK X, \cK Y)\; : 
         \qquad U \mapsto \bigcup_{g \in G} g(U) \, .
\end{equation}
\begin{prop} \label{unionCompactLip}
    ${}_{}^\cup\! \cK G \in C(\cK X, \cK Y)$ is Lipschitz with Lipschitz 
    constant 
    $$\hphantom{xxxxxxxxxxx}
    r^{}_{{}_{}^\cup\! \cK G} \, \le \, \sup \{r^{}_{\cK (g)} \mid g\in G\}
    \, = \, \sup \{r^{}_g \mid g\in G\} \, \le \, r \, . 
    \hphantom{xxxxxxxxx}\qed $$ 
\end{prop} 

If $G$ consists of contractions and $X = Y$ then we call $G$ a {\em compact 
iterated function system} (IFS). \index{compact IFS}\index{IFS}
This generalizes the previous definition 
which will now be referred to as a {\em finite} IFS. Hutchinson's
theorem evidently generalizes:
\begin{prop} \label{family-contra}
   Let  $G$ be a compact admissible family of contractions from $X$ to $Y$ 
   with uniform Lipschitz bound $r <1$. Then 
   ${}_{}^\cup\! \cK G\!: \cK X \longrightarrow \cK Y$ is a contraction with 
   Lipschitz constant at most $r$. If, in addition, $X=Y$
   $\,(\/$so $G$ is a compact {\rm IFS}$)$, then there is a unique
   $W \in \cK X$ $\,(\/$the attractor of the {\rm IFS}$)$ which is invariant 
   under $G$\! :  $W = \bigcup_{g \in G} \,  g(W)$. For arbitrary $Z \in \cK X$, 
   the iterates $({}_{}^\cup\!\cK G)^{\ell}(Z)$ converge to $W\in\cK X$
   as $\ell\to\infty$.
\end{prop}
{\sc Proof:} For  $Z \in \cK X$, 
$\bigl(({}_{}^\cup\!\cK G)^{\ell}(Z)\bigr)_{\ell\in\NN}$ is a Cauchy 
sequence in $\cK X$ since ${}_{}^\cup\! \cK G$ is a contraction.
Completeness shows the existence of a limit, which is a fixed point, and 
uniqueness follows immediately from contractivity.  \qed

\section{Self-similar measures for compact sets of contractions}
\label{sec-3}

One of the important contributions of Hutchinson \cite{Hutchinson} was the 
realization that the attractor of an IFS carries a {\em measure} that 
\index{measure} is likewise 
an invariant of the IFS, and indeed a far finer one than the attractor itself. 
In this Section, we re-establish this result in the context of compact iterated 
function systems.  The basic assumption that we need is that our compact space 
of contractions carries a measure of its own. We do not need to specify in
advance what this measure is. In the original case studied by Hutchinson,
where the IFS was finite, this supplementary measure was (effectively)
counting measure.

Let $X$ be a compact metric space. We denote by $\cP (X)$ the space of all 
probability measures on $X$ --- that is, positive regular Borel measures $\mu$ 
with total measure $\mu(X)=1$, see \cite[Ch.\ IV.4]{RS} for background 
material. Note that 
$$\cP(X)\;\subset\;\cM_+(X)\;\subset\;\cM(X)\;\subset\;\cM^{}_{\CC}(X)\, , $$ 
where $\cM_+(X)$, $\cM(X)$, and $\cM^{}_{\CC}(X)$ denote the spaces of
positive, signed (or real), and complex regular Borel measures, respectively.
For later use, we also define
$$\cM_+^m(X) \; := \; \{\mu\in\cM_+(X)\mid\mu(X)=m\}$$
so that $\cP(X)=\cM_+^1(X)$.
Since the Riesz-Markov theorem \cite[Thms.\ IV.14 and IV.18]{RS} states
that regular Borel measures are in one-to-one correspondence with linear 
functionals on the space $C(X,\RR)$ (resp.\ $C(X,\CC)$), equipped
with the compact-open topology, we shall usually
identify these pictures. So, we shall write $\mu(E)$ for the measure of a 
Borel set $E$, but often also $\mu(g)$ instead of $\int_X g \,{\rm d}\mu$
for the measure ($=$ integral) of a function.

In view of this, it is natural to
equip $\cM(X)$ (resp.\ $\cM^{}_{\CC}(X)$) with the weak-$*$ topology,
called the {\em vague topology\/} in this context \cite[p.\ 114]{RS},
\index{vague topology} the weakest topology that makes all the mappings 
$\mu \mapsto \mu(g)$ of $\cM (X) \to \RR$ 
(resp.\ of $\cM^{}_{\CC} (X) \to \CC$) continuous, where $g\in C(X,\RR)$
(resp.\ $g\in C(X,\CC)$). Let us mention that, since $X$ is compact, 
$C(X,\CC)$ is actually a Banach space with the sup-norm $\|.\|_{\infty}$
(which induces the compact-open topology), hence $\cM^{}_{\CC}(X) = C(X,\CC)^*$ 
can also be viewed as a Banach space, with induced norm 
$$\|\mu\|\; := \;\sup\{\, |\mu(g)|\,\mid g\in C(X,\CC),\, \|g\|_{\infty}=1\}\,$$
(see \cite[Thm.\ III.2]{RS}). 
With Hahn's decomposition theorem \cite[Thm.\ IV.16]{RS}, one
gets $\|\mu\|=|\mu|(X)$, where $|\mu|\in\cM_+(X)$ 
denotes the total variation measure of $\mu$. The analogous statement also 
holds for $\cM(X)=C(X,\RR)^*$. We shall need both ways to look at $\cM(X)$.
Another result in this context, using the Banach-Alaoglu Theorem
\cite[Thm.\ IV.21]{RS}, is that the unit balls in $\cM^{}_{\CC}(X)$ and
$\cM(X)$ are compact in the vague topology. The closed subspace
$\cP(X)=\{\mu\in\cM_+(X)\mid\|\mu\|=1\}$ is then also compact.  
{}For an alternative derivation of the last statement, without reference 
to the Banach-Alaoglu Theorem, see \cite[Prop.~8.1]{Akin}.

{}Following Hutchinson, we now define a {\em metric} 
\index{Hutchinson metric} on $\cP(X)$:
\begin{equation}\label{metric-def}
   L(\mu,\nu) \; := \; \sup \{\, |\mu(\phi) - \nu(\phi)| \,
             \mid \phi \in {\rm Lip}(1,X,\RR) \}\, .
\end{equation}
In fact, it is hardly clear that this {\em is} a metric\footnote{
This is well-known among experts, see \cite[{\S} 4.3]{Hutchinson}, 
but we could not find an explicit proof in the literature. 
Since it is an important part of our argument and the proof is not 
entirely trivial, we include it here. Note that the restriction to 
$\cP(X)$ (or to $\cM_+^m(X)$ for some $m>0$) is vital.}, 
but we shall show that below. 
It is useful to observe that, in this definition, we can also replace 
${\rm Lip}(1,X,\RR)$ by ${\rm Lip}(\le\! 1,X,\RR)$ without altering the 
resulting function $L$.  We often make use of this in the sequel.
Since $X$ is compact, ${\rm diam}(X):=\sup\{d(x,y)\mid x,y\in X\}$,
the {\em diameter} \index{diameter} of $X$, is finite, and we can
state another property explicitly.
\begin{lemma}\label{alt-metric}
   Let $L$ be defined, on $\cP(X)$, by Eq.~{\rm (\ref{metric-def})}. 
   Then we have
   $$L(\mu,\nu) \; = \; \sup \{\, |\mu(\psi) - \nu(\psi)| \,
   \mid \psi \in \cL \}$$ where $\cL:=\{\psi \in {\rm Lip}(1,X,\RR) \mid 
   ||\psi||_{\infty} \le {\rm diam}(X) \}$. 
   Furthermore, $\cL \subset C(X,\RR)$ is 
%  equicontinuous, and hence 
   compact $\,(\/$in the compact-open topology $\!)$.
\end{lemma}
{\sc Proof}:
Let $\mu,\nu \in \cP(X), \; \phi \in {\rm Lip}(1,X,\RR)$. 
Then, for any  $c\in \RR$, we have
$|\mu(\phi - c) - \nu(\phi -c) | = |\mu(\phi) - \nu(\phi)|$
since $\mu(c) = \nu(c) \, (=c)$. Let $a \in X$. Choosing $c=\phi(a)$ we obtain 
$|\phi(x) -\phi(a)| \le r^{}_\phi d(x,a) \le {\rm diam}(X)$.
So, $\psi(x):=\phi(x) - \phi(a)$ is a function in $\cL$, and the 
restriction to $\cL$ does not change the supremum value of 
$|\mu(\phi)-\nu(\phi)|$.  This establishes the first assertion.

Note that $\cL$ is clearly closed in $C(X,\RR)$.
Also, since the $\psi\in\cL$ are uniformly bounded, we have 
$\overline{\cL(x)} = \overline{\{\psi(x)\mid\psi\in\cL\}} 
\subset [-{\rm diam}(X), {\rm diam}(X)]$
for every $x\in X$, so each $\overline{\cL(x)}$ is compact. 
Finally, $\cL$ is equi-continuous since it consists of Lipschitz functions 
with uniformly bounded Lipschitz constants. By Ascoli's theorem, 
see \cite[Thms.\ 7.21 and 7.22]{Kelley}, $\cL$ itself is then
compact in $C(X,\RR)$ (in the compact-open topology). \qed

\begin{prop} \label{metricL}
  $L$ is a metric on $\cP(X)$ and induces the vague topology on $\cP(X)$.
  In particular, $\cP(X)$ is a complete metric space.
\end{prop}

{\sc Proof}:
That $\cP(X)$ is a complete metric space follows from its compactness
(see above) as soon as we have shown $L$ to be a metric.  

We use Lemma~\ref{alt-metric}.
If $\mu,\nu \in \cP(X)$, then $\mu-\nu : \cL \longrightarrow \RR$
is continuous and so has compact image. This shows that $L(\mu,\nu)$
is finite. Non-negativity and symmetry are obvious,
as is the triangle inequality. Thus $L$ is certainly a pseudo-metric. 
It remains to be shown that $L(\mu,\nu)=0$ implies $\mu=\nu$. 
Assume the converse and set $\omega=\mu-\nu$. Then
$c:=\|\omega\|>0$ and there is a $\delta>0$ and a function 
$g\in C(X,\RR)$ such that $|\omega(g)|\ge\delta>0$.

Since Lipschitz functions 
are dense in $C(X,\RR)$, see Lemma~\ref{lip-real} in the Appendix, we can 
choose $\phi$ Lipschitz with $\|g-\phi\|_{\infty}<\delta/2c$. Then, 
$$|\omega(g-\phi)| \; \le \; \|\omega\| \cdot \|g-\phi\|_{\infty}
   \; < \; \frac{\delta}{2} $$
and thus $|\omega(\phi)|\ge \delta/2 > 0$. Now, we don't know the
Lipschitz constant $r^{}_{\phi}$ of $\phi$, but
$\phi':=\phi/r^{}_{\phi}$ is in ${\rm Lip}(1,X,\RR)$ and still
$|\omega(\phi')|>0$, so $L(\omega,0)=L(\mu,\nu)>0$, which contradicts
the assumption. This shows that $L$ is a metric.

{}Finally, we compare the topologies. Let $\mu_n \rightarrow \mu$ vaguely 
as $n\to\infty$. Since $\mu_n\in\cP(X)$, we get 
$|\mu_n(\phi_1)-\mu_n(\phi_2)|\le\|\phi_1 -\phi_2\|_{\infty}$
independently of $n$, so the $\mu_n$ constitute a family of equi-continuous
mappings in $C(C(X,\RR),\RR)$. Consequently, by Lemma~\ref{point-compact} of
the Appendix, $\mu_n(\phi)\rightarrow\mu(\phi)$ uniformly on $\cL$ because 
$\cL$ is compact (Lemma~\ref{alt-metric}), and hence $L(\mu_n,\mu)\to 0$.
% Conversely, suppose that $L(\mu_n,\mu)\to 0$ as $n\to\infty$. 
% {}From the compactness
% of $\cP(X)$ in the vague topology, $\{\mu_n\}$ has a convergent subsequence, 
% $\{\mu_{n_i}\}$ say, with limit $\omega \in \cP(X)$. But then 
% $L(\mu_{n_i},\omega)\to 0$ by the first part, 
% so $\omega = \mu$. Since this is true
% for {\em every} convergent subsequence of $\{\mu_n\}$, we see that 
% $\mu_n\to \mu$ also in the vague topology. 
Conversely, observe that $\cL$ is compact in the vague topology and
Hausdorff under the metric $L$. But the identity is a one-to-one mapping,
and (due to the previous argument) also continuous, when
viewed as a mapping from $\cL$ with the vague topology to $\cL$ with the
metric topology. Therefore, it is a homeomorphism \cite[Thm.\ 5.8]{Kelley}, 
and the topologies coincide.       \qed

\smallskip
Let $F$ be a compact IFS on the compact
metric space $X$ and let $W \in \cK X $ be its attractor. Let
$r^{}_F := \sup\{r^{}_f \mid f \in F \} < 1$.

{}For each $f\in F$, $f(W) \subset W$, and we obtain a bounded linear operator
$f.(\centerdot)$ on the space of all signed Borel measures $\mathcal{M}(W)$ 
of $W$  by
\begin{equation} \label{measureMap}
   \mu \mapsto f.\mu \; ; \quad f.\mu(\phi) \; := \; \mu(\phi \circ f)
\end{equation}
for all $\phi \in C(W,\RR)$.
Evidently, $f.(\centerdot):\; \cP(W) \longrightarrow \cP(W)$, i.e.\
if $\mu$ is a probability measure, so is $f.\mu$. Note that the matching
definition for Borel sets $E$ reads $f.\mu(E) = \mu(f^{-1}(E))$ where
$f^{-1}(E)$ is the preimage of $E$ under $f$.

\begin{prop} 
   The mapping $F \times \cP(W) \longrightarrow  \cP(W)$
   defined by $(f,\mu) \mapsto f.\mu$ is continuous.
\end{prop}
{\sc Proof}: $\cP(W)$ is a compact metric space, so certainly Hausdorff.
Fix $f\in F$ and consider the mapping $\mu\mapsto f.\mu$. Then, for
$\phi\in {\rm Lip}(1,W,\RR)$, we clearly have
$|f.\mu(\phi)-f.\nu(\phi)|=|(\mu-\nu)(\phi\circ f)|
\le L(\mu,\nu)$ because $\phi\circ f$ has Lipschitz constant $\le 1$
due to the definition of $F$. So $\mu\mapsto f.\mu$ is Lipschitz
and thus (uniformly) continuous on $\cP(W)$.
It follows that $F \times \cP(W) \longrightarrow \cP(W)$ is 
jointly continuous \cite[Thm.\ 7.5]{Kelley}. \qed 

\begin{prop} \label{contraction}
  Let $F$ be a compact IFS on the compact metric space $X$, with 
  attractor $W\in\cK X$. Let $\nu \in \cP(F)$. Then the
  $\nu$-averaging  mapping 
\begin{equation}
   \begin{array}{rccl}
   \mathcal{A}_\nu: & \cP(W) & \longrightarrow & \cP(W) \\
                    & \mu    & \mapsto         & 
                      \int_F (f.\mu) \, {\rm d}\nu(f)
   \end{array}
\end{equation}
   is a contraction relative to the $L$-metric on $\cP(W)$, 
   with contraction constant at most 
   $r:= r^{}_F=\sup\{r^{}_f\mid f\in F\}$. 
   In particular, there is a unique measure 
   $\rho^{(\nu)} \in \cP(W)$ which satisfies 
\begin{equation}
   \rho^{(\nu)} \; = \; \int_F (f.\rho^{(\nu)}) \, {\rm d}\nu(f).
\end{equation}
\end{prop}
{\sc Proof}: Let $\omega_1,\omega_2 \in \cP(W)$ and let 
$\phi\in{\rm Lip}(1,W,\RR)$. Then 
\begin{eqnarray*}
  \bigl|(\mathcal{A}_\nu (\omega_1))(\phi)  
        - (\mathcal{A}_\nu (\omega_2))(\phi)\bigr|
  & = & \left|\int_F\omega_1(\phi\circ f) \,{\rm d}\nu(f) 
       - \int_F\omega_2(\phi\circ f) \,{\rm d}\nu(f) \right| \\
  & = & \left|\int_F \,\Bigl(\omega_1(\phi\circ f) 
       - \omega_2(\phi\circ f) \Bigr)\,{\rm d}\nu(f) \right| \\
  &\le& r \int_F  \left|\omega_1(r^{-1}\phi\circ f) 
      - \omega_2(r^{-1}\phi\circ f)  \right|\,{\rm d}\nu(f)  \\
  &\le& r \int_F L(\omega_1,\omega_2) \, {\rm d}\nu(f) 
          \;\, = \;\, r \, L(\omega_1,\omega_2) \; \nu(F) \\
  & = & r \, L(\omega_1,\omega_2)\, ,
\end{eqnarray*}
since $r^{-1}\phi\circ f \in {\rm Lip}(\le 1,W,\RR)$.
This being true for all $\phi \in {\rm Lip}(1,W,\RR)$, we have 
\begin{equation}
   L(\mathcal{A}_\nu (\omega_1), \mathcal{A}_\nu (\omega_2)) 
   \; \le \;  r \, L(\omega_1,\omega_2)
\end{equation}
which is what we wanted. The existence of the unique measure
$\rho^{(\nu)}$ follows now, once again, from the general contraction
principle. \qed

\smallskip

{\sc Remark}: More generally, we may replace $W$ in Proposition \ref{contraction}
by any other $W^+\in\cK X$ that satisfies $f(W^+)\subset W^+$ for all $f\in F$.
Of course, $W\subset W^+$ and the invariant measure for $W^+$
is supported on $W$, hence is effectively the same as $\rho^{(\nu)}$.

\section{Affine mappings in locally compact Abelian groups}\label{affine-maps}

In this Section, we treat the foregoing material in the setting of a
locally compact Abelian group \index{Abelian group!locally compact}
(LCAG). A convenient source for the results on LCAGs that
we need is \cite[Ch.\ 1]{Rudin}.

\subsection{Affine mappings}
Let $H$ be an additive LCAG \index{LCAG} that is equipped
with a translation invariant metric $d$ with respect to which $H$ is complete.
For more information on metrizability, see \cite[Ch.\ 2, {\S} 8]{HR}.
We also assume that an automorphism $A$\index{automorphism} is given (in particular,
$A(0)=0$) and that a Haar measure $\theta$ on $H$ has been fixed. 
\index{Haar measure} It is unique up to normalization.

Since $A.\theta := \theta\circ A^{-1}$ is also $H$-invariant, it is another 
Haar measure, and we thus have
$A.\theta = \alpha \theta$ for some $\alpha > 0$, the {\em modulus} of $A$.
\index{modulus} If $H=\RR^n$, $\theta$ is Lebesgue measure, $A$ is simply an 
invertible linear map, and $\alpha=|\det(A^{-1})|=1/|\det(A)|$. 
\begin{lemma}\label{contraction-modulus}
   Let $H$ be as described with metric $d$.
   If $A$ is a contraction on $H$ relative to $d$, then $A$ has
   modulus $\alpha > 1$.
\end{lemma}
{\sc Proof}: Since $H$ is locally compact, it contains a compact neighbourhood
$U$ of $0$. By assumption, the topology given on $H$ agrees with the metric
topology. So, we can choose $U$ such that inside it we can find balls $B_r(0)$ 
and $B_s(0)$ with $r>s>0$ and the property that $B_r(0)\setminus B_s(0)$
contains a non-empty open set which must then have positive measure.

On the other hand, $A$ is a contraction, $A(0)=0$, and 
$d(0,A^n(U))\to 0$ as $n\to\infty$ for any compact $U\subset H$.
So, there must be some $m\in\NN$ such that $A^m(B_r(0))\subset B_s(0)$. Combined
with the previous argument, this says $A^m(B_r(0))$ has smaller measure than 
$B_r(0)$. Consequently, the modulus of $A^m$ is $\alpha^m > 1$, and then
also $\alpha > 1$. \qed

\smallskip
Clearly, the converse of Lemma \ref{contraction-modulus} is not true.
In view of the general context of this paper, we assume from now on
that $A$ is a contraction on $H$ relative to $d$.
% \begin{itemize}
% \item[\bf A1] $A$ is a contraction on $H$ relative to $d$.
% \item[\bf A2] The modulus of $A$ fulfils $\alpha > 1$.
% \end{itemize}
A mapping $f:\; H\to H$ of the form
\begin{equation}
   f: \; x \mapsto A(x) + v
\end{equation}
with $v\in H$ is called an {\em affine map} \index{affine map}
with {\em automorphism} $A$ \index{automorphism} and {\em translation}
$v$. \index{translation} We will sometimes denote this mapping by $A_v$.

Let $\cM^{}_{\CC}(H)$ be the space of all bounded regular complex Borel 
measures $\lambda$ on $H$, i.e.\ measures with 
$\|\lambda\|=|\lambda|(H) < \infty$.  We recall that the
convolution of two measures $\lambda_1,\lambda_2 \in \cM^{}_{\CC}(H)$
is defined by
\begin{equation}
     (\lambda_1 * \lambda_2)(\phi) \; = \;
     \int_{H\times H} \phi(x+y) 
     \, {\rm d}\lambda_1(x) \, {\rm d}\lambda_2(y) \, ,
\end{equation}
for $\phi\in C(H,\CC)$. If formulated for Borel sets $E$, the matching 
equation is
\begin{equation}
     (\lambda_1 * \lambda_2)(E) \; = \;
     (\lambda_1 \otimes \lambda_2)(E^{(2)})
\end{equation}
where $E^{(2)} :=  \{ (x,y) \in H\times H\mid x+y \in E \}$.

The Fourier-Stieltjes transform of $\lambda\in\cM^{}_{\CC}(H)$ is the
function $\widehat{\lambda}$ defined on the dual group $\widehat{H}$ of $H$ by
\begin{equation}\label{fTmeas}
   \widehat{\lambda}(k) \; = \; 
        \int_H \overline{\langle k,x \rangle} \, {\rm d}\lambda(x)
\end{equation}
where $x\mapsto\langle k,x\rangle$ is the continuous {\em character} on $H$ 
defined by $k\in\widehat{H}$. \index{character}

The automorphism $A$ on $H$ induces an automorphism $A^T$ on $\widehat{H}$:
$k\mapsto A^T k$ where $A^T k$, in turn, defines the character 
$x \mapsto \langle k,Ax\rangle$ on $H$.

We collect now some basic facts that we need. These are all elementary
consequences of the definitions, whence we omit proofs. We write $A.h$ for 
the function defined
by $x\mapsto h(A^{-1}(x))$ in analogy to $A.\mu=\mu\circ A^{-1}$ for 
measures, and $h\mu$, with $h\in L^1(H)$, for the measure defined by
$(h\mu)(\phi)=\mu(h\phi)$. Thus, $h\mu$ is absolutely continuous 
with respect to $\mu$, and $h$ is the corresponding Radon-Nikodym 
density (also called Radon-Nikodym derivative).
\begin{prop} \label{facts}
  Let $H,\theta,A,\alpha$ be as defined above. Let 
  $\lambda_1,\lambda_2\in\cM^{}_{\CC}(H)$ and $h\in L^1(H)$.
  Then we have
\begin{enumerate}
\item $ {\rm d}\theta(A^{-1}x) \, = \, \alpha \, {\rm d}\theta(x)$
\item $ A.(h \theta) \, = \, \alpha (A. h)\theta$
\item $ A.(\lambda_1 * \lambda_2) \, = \, A.\lambda_1 * A.\lambda_2 $
\item $ \widehat{A.\lambda} \, = \, \widehat{\lambda}\circ A^T 
        \, = \, (A^T)^{-1}.\widehat{\lambda}$     \qed
\end{enumerate}
\end{prop}

\subsection{Compact families of affine mappings} \label{cfam}
Assume now that $F$ is a compact family of contractions on the
LCAG $H$, each $f\in F$
being of the form $$ A_v : \; x\mapsto Ax+v $$ for some $v\in H$ (but all
having the same $A$, namely our contractive automorphism fixed above).
Evidently, $F$ is a compact admissible family of mappings from $H$ to $H$.
Define
\begin{equation}
    F_H \; = \; \{ v\mid A_v\in F\}
        \; = \; \{ f(0)\mid f\in F\} \; = \; F(0)
        \; \subset \; H \, .
\end{equation}
In view of (\ref{evaluation}) and (\ref{extension}), the mapping $F\to F_H$ 
induced by $f\mapsto f(0)$ is continuous and hence $F_H$ is compact and 
homeomorphic to $F$. In particular, there is a natural isomorphism between
$\cM^{}_{\CC}(F)$ and the space of regular measures on $H$ that are
supported on $F_H$.

Let $\nu^{}_F\in\cP(F)$ and $\nu\in\cP(H)$ be such a corresponding pair of
(probability) measures. We then have an averaging operator
$\cA_{\nu}:\; \cP(W^+)\to\cP(W^+)$ whenever $W^+\subset H$ is any compact
subset of $H$ for which $FW^+\subset W^+$.

Let $\lambda\in\cP(W^+)$. Then, for all Borel sets $E\subset H$,
\begin{eqnarray}
  \cA_{\nu} \lambda(E) & = &
      \int_F f.\lambda(E)\, {\rm d}\nu^{}_F(f)
      \; = \; \int_F \lambda(f^{-1}(E))\, {\rm d}\nu^{}_F(f)
      \nonumber \\
       & = & \int_H \lambda(A^{-1}(E-v))\, {\rm d}\nu(v) 
       \; = \; \int_H (A.\lambda)(E-v)\, {\rm d}\nu(v)
       \nonumber \\
       & = & \int_{H\times H} \bs{1}^{}_E(u+v) 
             \, {\rm d}(A.\lambda)(u) \, {\rm d}\nu(v) \nonumber \\
       & = & \vphantom{\int}(\nu * A.\lambda)(E)  \label{iter}
\end{eqnarray}
where $\bs{1}^{}_E$ is the characteristic function of $E$.
Since $\textrm{supp} (\lambda) \subset W^+$, we have
$f(\textrm{supp} (\lambda)) \subset W^+$ for all $f\in F$.
This implies $\textrm{supp} (\cA_{\nu} \lambda) \subset W^+$
and, more generally, $\textrm{supp} (\cA_{\nu}^{\ell} \lambda) \subset W^+$
for all $\ell\ge 0$. In particular, we can also infer
$$\textrm{supp} (\cA_{\nu} \lambda) \; = \;
  \textrm{supp} (\nu * A.\lambda)   \; \subset \;
  \textrm{supp} (\nu) + A \, \textrm{supp} (\lambda)
  \; \subset \; W^+ \, . $$
It is clear that we can now iterate (\ref{iter}) to get
$$ \cA_{\nu}^{\ell} \lambda \; = \;
   \nu * A.\nu * \ldots * A^{\ell-1} .\nu * A^{\ell} .\lambda $$
together with the inclusion relation
$$\textrm{supp} (\cA_{\nu}^{\ell} \lambda) \; \subset \;
  \textrm{supp}(\nu) + A \, \textrm{supp}(\nu) + \ldots +
  A^{\ell-1} \, \textrm{supp}(\nu) + A^{\ell} \, \textrm{supp}(\lambda)
  \; \subset \; W^+ \, . $$
Since $A$ is a contraction, 
$A^{\ell} \, \textrm{supp}(\lambda)\longrightarrow\{0\}$ (in $\cK H$)
as $\ell\to\infty$ and we have
$$ \sum_{\ell=0}^{\infty} A^{\ell} \, \textrm{supp}(\nu)
   \; \subset \; W^+ \, .$$
In particular, since $W$ is $F$-invariant, we must also have
$$ \sum_{\ell=0}^{\infty} A^{\ell} \, \textrm{supp}(\nu)
   \; \subset \; W \, .$$
Define $W^+$ to be the smallest compact subset of $H$ which is $F$-invariant
and contains $\sum_{j=0}^{\ell} A^j\,\textrm{supp}(\nu)$ for all $\ell\geq 0$.

Define $\omega^{(\ell)}\in\cP(W^+)$ by $\omega^{(0)}=\nu$, and (for $\ell\ge 0$)
$$\omega^{(\ell+1)} \; = \; \cA_{\nu} \omega^{(\ell)}
  \; = \; \nu * A.\omega^{(\ell)}\, .$$
Next, let $\omega$ be the unique $\cA_{\nu}$-invariant measure of
$\cP(W^+)$, see Proposition \ref{contraction} and the Remark following it. 
We know, again by Proposition \ref{contraction}, that $\cA_{\nu}$
is a contraction on $\cP(W^+)$. Moreover, by (\ref{iter}), 
$\{\omega^{(\ell)}\}$ contracts to $\omega$ as $\ell\to\infty$. Thus
$$\lim_{\ell\to\infty} \omega^{(\ell)} \; = \; \omega$$ 
in $\cP(W^+)$, with convergence in the vague topology.

\begin{prop} Under the above assumptions, we have
\begin{enumerate} \label{fix-measure}
\item $\omega = \conv_{\ell=0}^{\infty} (A^{\ell} .\nu)\in\cP(W)$,
      which converges in the vague topology, 
      is the unique self-similar probability measure for the compact 
     admissible family of contractions $F=\{A_v\mid v\in F_H\}$ 
      with respect to
      the measure $\nu^{}_F$ on $F$.\index{measure!self-similar}
\item $\widehat{\omega} = \prod_{\ell=0}^{\infty} (A^T)^{-\ell} .
      \widehat{\nu}$, convergence being uniform convergence on compact
      sets $\,(\/$compact convergence$)$. \index{convergence!compact}
\item If the convolution product for $\omega$ converges also in the
      $\|.\|$-topology on $\cP(W)$, the convergence of $\widehat{\omega}$ is
      $\,(\/$globally $\!)$ uniform. \index{convergence!uniform}
\end{enumerate}
\end{prop}
{\sc Proof}: Part 1 is clear from the discussion above. The support of
$\omega$ is inside $W$ by Proposition \ref{contraction}. Part 2 follows from 
Proposition \ref{facts} and the continuity of the Fourier transform,
sending measures $\mu$ to bounded and uniformly continuous functions
$\widehat{\mu}$. The convergence statement is a direct consequence of
L\'evy's continuity theorem, see Theorem~\ref{levy} of the Appendix.
Finally, the third assertion follows directly from 
$\|\widehat{\mu}\|^{}_{\infty}\le\|\mu\|$, see (\ref{fTmeas}),
without reference to Part 2. \qed

\smallskip
{\sc Remark}: It is only a matter of convenience to start the above iteration
with $\omega^{(0)}=\nu$. Any other choice $\lambda\in\cP(W^+)$ is 
equally admissible and will lead to the same result, because 
$A^{\ell}.\lambda\to\delta_0$, as $\ell\to\infty$, and $\delta_0$,
the unit point measure at $0$, is the neutral element of convolution,
i.e.\ $\mu*\delta_0=\mu$ for all measures $\mu$.

\subsection{Self-similar functions} \label{self-sim-functions}
If we assume that the measure $\nu$ on our compact family of affine
contractions is absolutely continuous with respect to Haar measure, then
Proposition \ref{fix-measure} gets re-interpreted in terms of functions rather
than measures.

We suppose the same notation as in Section \ref{cfam} and assume in
addition that the measure $\nu$ derived from $\nu^{}_F$ on $F$ is 
absolutely continuous w.r.t.\ $\theta$, so $\nu$ is of the 
form $\nu = h\theta$, where $h\in L^1(H)$ and $\textrm{supp}(h)\subset F_H$,
with $F_H$ compact. For such measures, convolution
matches the usual convolution of functions. Thus, using Proposition \ref{facts}
and Part 1 of Proposition \ref{fix-measure}, we obtain
\begin{eqnarray*}
    \omega^{(\ell)} & = &
     \nu * A.\nu * \ldots * A^{\ell}.\nu \\
    & = & h\theta * \alpha (A.h)\theta * \ldots *
          \alpha^{\ell} (A^{\ell}.h)\theta \\
    & = & \left( \conv_{j=0}^{\ell} \, \alpha^j (A^j .h)\right) \theta
\end{eqnarray*}
and $\omega = \bigl( \conv_{j=0}^{\infty}\,\alpha^j (A^j .h)\bigr) \theta$,
with convergence so far only in the vague topology. 
However, as the brackets already imply, convergence in the
$\|.\|$-topology would be preferable.  The situation is as follows.
If we identify $L^1(H)$ with a subspace of $\cM^{}_{\CC}(H)$ via \index{$L^1$}
$f\mapsto f\theta$, this is a {\em closed subspace} of $\cM^{}_{\CC}(H)$
in the $\|.\|$-topology, see \cite[{\S}  1.3.5]{Rudin}. Consequently, the
$\|.\|$-convergence of absolutely continuous measures is equivalent
to the $L^1$-convergence of their Radon-Nikodym densities in $L^1(H)$.

To establish also the $\|.\|$-convergence in our case, recall first the
following result \cite[Thm.\ 1.1.8]{Rudin} on approximate units in
the commutative convolution Banach algebra $L^1(H)$ (with norm $\|.\|^{}_1$).
\index{convolution algebra} \index{approximate unit}
\begin{lemma} \label{approx-unit} 
  Given $f \in L^1(H)$ and $\varepsilon > 0$, there exists a neighbourhood
  $V$ of $0$ in $H$ with the following property: if $u$ is a non-negative
  Borel function which vanishes outside $V$, and if
  $\,\| u \|^{}_1 = \int_H u(x)\, {\rm d}\theta(x) = 1$, then
  $$ \hphantom{xxxxxxxxxxxxxxxxxxxxxxxx}
  \|\, f - f * u \,\|^{}_1 \; < \; \varepsilon \, . 
     \hphantom{xxxxxxxxxxxxxxxxxxxxxx}\qed$$
\end{lemma}

We are now in the following situation. Our starting function is
$h\in L^1(H)$, with ${\rm supp}(h)$ compact, $h\ge 0$ and
$\int_H h \, {\rm d}\theta = \|h\|^{}_1 = 1$. Let
$f^{}_{\ell} = \alpha^{\ell} (A^{\ell}.h)$ for $\ell\ge 0$, so that
$f^{}_{\ell}\ge 0$ and $\|f^{}_{\ell}\|^{}_1=1$. Also, 
${\rm supp}(f^{}_{\ell}) = A^{\ell} \, {\rm supp}(h)$, and
we have the relation
$\|f^{}_{\ell} * f^{}_{\ell+1} * \ldots * f^{}_{\ell+k}\|^{}_1 = 1$
for all $k\ge 0$.
\begin{prop} \label{norm-conv} \index{norm topology}
  Let $F$ be a compact family of affine mappings, with contractive automorphism
  $A$, modulus $\alpha$ and attractor $W\subset H$. Let $F_H$ be the
  corresponding set of translations and let $\nu=h\theta$ be an absolutely
  continuous probability measure on $F_H$, where $h\in L^1(H)$.
  Then, the infinite convolution product
  $\conv_{\ell=0}^{\infty} \, f^{}_{\ell}$ converges to an $L^1$-function,
  hence $\conv_{\ell=0}^{\infty} \, f^{}_{\ell} \theta$ converges
  also in the $\|.\|$-topology.
\end{prop}

{\sc Proof}: Since $L^1(H)$ is complete, it suffices to show that
$\bigl(\conv_{\ell=0}^{n} \, f^{}_{\ell}\bigr)_{n\ge 0}$ is 
a Cauchy sequence.
Fix $\varepsilon>0$ and let $V$ be the neighbourhood for 
the $L^1$-function $f=f^{}_0=h$
according to Lemma \ref{approx-unit}. Since $A$ is a contraction, there
exists an integer $N$ so that 
$\sum_{\ell\ge N} {\rm supp}(f^{}_{\ell}) \subset V$.
In particular, any finite convolution of the form
$\conv_{\ell=N}^{N+k} \, f^{}_{\ell}$, $k\ge 0$,
is then an approximate unit for $h$ with bound $\varepsilon$.

Let now $n,m \ge N$ and define $u=\conv_{\ell=N}^{n} \, f^{}_{\ell}$
and $v=\conv_{\ell=N}^{m} \, f^{}_{\ell}$. Then
\begin{equation}\nonumber
\begin{split}
    \left\|\,
    \overset{n}{\underset{\ell=0}{\conv}} \, f^{}_{\ell} - 
    \overset{m}{\underset{\ell=0}{\conv}} \, f^{}_{\ell} \, \right\|_1 
  & \; = \;
    \left\| \left( 
    \overset{N-1}{\underset{\ell=0}{\conv}}\, f^{}_{\ell}\right)
    * (u-v) \right\|_1 \\
  & \; \le  \;
    \left\| \,
    \overset{N-1}{\underset{\ell=1}{\conv}}\, f^{}_{\ell} \, \right\|_1
    \cdot \bigl\| \, h * u - h * v \, \bigr\|_1   \\
  & \; = \;
    \bigl\| \, (h * u - h) + (h - h * v) \, \bigr\|_1 \\
  & \; \le \;
    \bigl\| \, h - h * u \, \bigr\|_1 + 
    \bigl\| \, h - h * v \, \bigr\|_1
    \; < \; 2 \varepsilon \\
\end{split}
\end{equation}
by application of Lemma \ref{approx-unit}. This gives part one of the claim,
while the rest follows, once again, from the Radon-Nikodym theorem.   \qed

\begin{prop}\label{single-conv}
  Under the general assumptions of Proposition {\em \ref{norm-conv}}, 
  we have:
\begin{enumerate}
\item There is a unique non-negative function $g\in L^1(H)$ which 
   satisfies\/\footnote{In \cite{JLS}, a mapping on functions of this form
   is called a continuous refinement operator. \index{refinement operator}}
   $$g \; = \; \alpha \int_H g(A^{-1}(x-v)) h(v)\, {\rm d}\theta (v)$$ 
   with normalization $\int_H g \, {\rm d}\theta = 1$.
\item $g = \conv_{\ell=0}^{\infty} \, \alpha^{\ell} (A^{\ell}.h)$,
   with convergence in the $L^1$-norm, and $\,{\rm supp}(g)\subset W$.
\item The Fourier transform of $g$ is the continuous function $\widehat{g} =
   \prod_{\ell=0}^{\infty} \alpha^{\ell} \mbox{\raisebox{0.2ex}{$\bigl($}}
   \widehat{h}.(A^T)^{\ell}\mbox{\raisebox{0.2ex}{$\bigr)$}}$,
   with uniform convergence of the product.
\item If $\vphantom{\prod_1^2} h\in L^1(H)\cap L^{\infty}(H)$, then $g$ is 
   continuous on $H$.
\end{enumerate}
\end{prop}
{\sc Proof}: 
The convergence claimed in Part 2 follows from Proposition \ref{norm-conv}, 
so $g\in L^1(H)$ and $\widehat{g}$ is then continuous.

{}From $\omega = \nu * A.\omega$, we then have, by Proposition \ref{facts},
$g = \alpha h * A.g$, which gives Part 1 by applying 
Proposition \ref{fix-measure}(1), and also the statement that
$\,{\rm supp}(g)\subset W$. 

The situation for Fourier transforms is even easier 
since $\widehat{g\theta} = \widehat{g}$ and we get the
product formula in Part 3 from Proposition \ref{fix-measure}(2) 
with uniform convergence by means of Proposition \ref{fix-measure}(3).

{}Finally, suppose that $h\in L^1(H)\cap L^{\infty}(H)$. 
Since $h\in L^{\infty}(H)$ and $A.g\in L^1(H)$, we obtain
(\cite[Thm.\ 1.1.6]{Rudin}) the continuity of $h*A.g$, hence
of $g$ itself.  \qed

\smallskip
{\sc Remark}: If $H=\RR^n$, we can actually iterate the last argument
and arrive at the stronger statement that 
$h\in L^1(\RR^n)\cap L^{\infty}(\RR^n)$ implies that $g$ is a
$C^{\infty}$-function with compact support contained in $W$. Furthermore,  if
$\nu$ is Lebesgue measure, then the
self-similar function $g$ enjoys remarkable properties with respect to
the averaging operator $\cA_\nu$, namely its partial derivatives are eigenfunctions for
the refinement operator with eigenvalues directly related to the
spectrum of $A$. This is the situation in our previously studied examples 
\cite{BM,BMmulti} and these results may be found there.

\section{Multi-component families of contractions}\label{mcms}

In this Section, we consider the generalization of the previous material to 
the case in which we have several compact metric spaces and sets of 
contractions between these spaces. This is the multi-component situation. 
The approach here is to consider the product of the various spaces in question.
The basic theorem on the existence of attractors then reduces at once to the
single-component situation already  dealt with. The question of self-similar 
measures also fits naturally into the product formalism, though the situation 
now acquires some new features that did not appear before.

\subsection{Contractions and attractors} 

Let $(X_1, d_1),  \dots, (X_n,d_n)$  be $n$ complete metric spaces and define
$N:=\{1,\dots,n\}$. Also, let $d^{}_{i,H}$ denote the corresponding
Hausdorff metric for $\cK X_i$, $i\in N$. We set
\begin{equation}
   X^{}_N \; := \; X_1 \times \ldots  \times X_n
\end{equation}
and write $x = (x_1, \dots, x_n)$ for the elements of $X^{}_N$.
We endow $X^{}_N$ with the metric
\begin{equation}
     d(x,y) \; := \; \sup \{d_i(x_i,y_i)\mid i\in N\}\, 
\end{equation}
relative to which it is also complete, and denote by $d^{}_H$ the attached
Hausdorff metric on $\cK X^{}_N$.

{}For each pair $(i,j)\in N\times N$, let $F_{ij}$  be a compact admissible 
family of contractions $f=f_{ij}:X_j \longrightarrow X_i$. 
We extend this to allow the 
possibility that $F_{ij}$ is empty, though we require that for each $i$ there 
is at least one $j$ for which $F_{ij} \ne \emptyset$. We let $0 < r < 1$ be a 
uniform upper bound on the contractivity factors of all these mappings. 
{}For each pair $(i,j)$, we have from (\ref{cKMap}) the mapping 
${}_{}^\cup\!\cK F_{ij} : \cK X_j \longrightarrow \cK X_i$. We define
\begin{equation}
   {}_{}^\cup\! \cK F : \quad \cK X_1 \times \ldots  \times \cK X_n \;
                \longrightarrow \; \cK X_1 \times \ldots  \times \cK X_n
\end{equation}
where
\begin{eqnarray} \label{multiContraction}
 {}_{}^\cup\! \cK F(U_1, \dots, U_n) 
   & :=& \Bigl(\mbox{$\bigcup_j$} \, {}_{}^\cup\! (\cK F_{1j}) \,(U_j), \dots ,
          \mbox{$\bigcup_j$} \, {}_{}^\cup\! (\cK F_{nj}) \, (U_j)\Bigr) 
          \nonumber \\
   &  =& \Bigl( \mbox{$\bigcup_j$} \mbox{$\bigcup_{f\in F_{1j}}$} 
                              f\, (U_j),  \dots ,
           \mbox{$\bigcup_j$} \mbox{$\bigcup_{f\in F_{nj}}$} 
                              f\, (U_j) \Bigr)\, .
\end{eqnarray}
Note that we write $(U_1, \dots, U_n)$ rather than $U_1\times\ldots\times U_n$
and that $\cK X_1 \times \ldots  \times \cK X_n$ is a strict subset of
$\cK X^{}_N$.

\begin{prop}\label{lip-contra}
  ${}_{}^\cup\! \cK F : \cK X_1 \times \ldots  \times \cK X_n \, 
  \longrightarrow \, \cK X_1 \times \ldots  \times \cK X_n$ is
  a contraction with Lipschitz constant at most $r$.
\end{prop}
{\sc Proof:}  For $U_i, V_i \in \cK X_i \,$, $i\in N$, we find
\begin{eqnarray*}
\lefteqn{ d^{}_H \,\bigl({}_{}^\cup\! \cK F(U_1, \dots, U_n), 
                    {}_{}^\cup\! \cK F(V_1, \dots, V_n) \bigr) } \\
&=& d^{}_H\,\bigl((\ldots,\,\mbox{$\bigcup_j$} \,{}_{}^\cup\! 
          (\cK F_{ij})\,(U_j), \ldots)\, , \;
         (\ldots,\,\mbox{$\bigcup_j$} \,{}_{}^\cup\! (\cK F_{ij}) \,(V_j),
          \ldots)\bigr) \\[1mm] 
&=& \sup_i\,\bigl\{d^{}_{i,H} \bigl(\mbox{$\bigcup_j$} \, 
          {}_{}^\cup\! (\cK F_{ij}) \,(U_j),
          \mbox{$\bigcup_j$} \, {}_{}^\cup\! (\cK F_{ij})\,(V_j)
          \bigr)\bigr\}  
    \qquad\qquad\quad \textrm{(by definition)}\\
&\le &  \sup_i \, \sup_j \,\bigl\{d^{}_{i,H}
        \bigl({}_{}^\cup\! (\cK F_{ij}) \,(U_j),
         {}_{}^\cup\! (\cK F_{ij}) \,(V_j)\bigr)\bigr\} 
    \quad\quad\qquad\qquad \textrm{(by (\ref{unionConst}))} \\ 
&\le& \sup_i \, \sup_j \,\{ r^{}_{\cK F_{ij}}\, d^{}_{j,H}(U_j,V_j) \}
      \; \le \; r \sup_j \, \{d^{}_{j,H}(U_j, V_j)\}
    \quad\quad \textrm{(by Prop.\ \ref{unionCompactLip})} \\
&=& r\,  d^{}_H \bigl((U_1, \dots , U_n), (V_1, \dots , V_n)\bigr) 
\end{eqnarray*} 
which establishes our assertion. \qed

\smallskip

We conclude, using the usual contraction principle, that there is a unique 
attractor for  ${}_{}^\cup\! \cK F$, in
$\cK X_1 \times \ldots  \times \cK X_n $, say $W_1 \times \ldots  \times W_n$. 
The $W_i$ thus 
form the unique solution (in compact sets) to the system of equations:
\begin{equation} \label{multiAttractor}
   W_i \; = \; \bigcup_{j=1}^{n} \; \bigcup_{f\in F_{ij}} f (W_j)
   \; , \quad i\in N \, .
\end{equation}

\subsection{Multi-component invariant measures}

The idea behind the invariant measures in the multi-component setting is 
straightforward in its conception but looks complicated in its details. 
We start with $n$ compact metric spaces $X_i$, $i\in N$, that are 
coupled by families $F_{ij}$ of contractions  
$f\!:X_j \longrightarrow X_i$. For the moment
we can take each set of mappings $F_{ij}$ to be finite. 

Each $f \in F_{ij}$ determines a transformation $\mu_j \mapsto f.\mu_j$ 
(see (\ref{measureMap}) for notation) of measure spaces
$\mathcal{M}(X_j) \longrightarrow \mathcal{M}(X_i)$.
Basically, we want to find a family of measures $\{\mu_1,\dots, \mu_n\}$ 
that is invariant under the average of these transformations:
\begin{equation} \label{basicInvariance}
   \mu_i \; = \; \sum_{j=1}^n \frac{1}{{\rm card} (F_{ij})} 
                 \sum_{f\in F_{ij}} f.\mu_j \, .
\end{equation}

There are some extensions and modifications that make this picture both 
more useful and easier to cope with mathematically:
\begin{enumerate}
  \item We are at liberty to give each set of mappings $F_{ij}$ its own
        weighting. 
  \item We need not restrict ourselves to {\em finite} sets $F_{ij}$, nor
     need we assume that our averaging is uniform within each of these sets.
     In what follows, we only assume that the sets $F_{ij}$ are compact spaces of
     mappings. We then deal with these points simultaneously by assigning
     positive\footnote{Strictly speaking, we should say non-negative measures,
     but we will always explicitly mention when the $0$-measure occurs.}
     measures $\sigma^{}_{ij}$ to each of these spaces of mappings. 
  \item It is mathematically easiest to deal with all of the measures
     $\{\mu_1, \dots , \mu_n\}$ as a single entity. Thus we prefer to deal with
     product measures $\mu_1 \otimes \ldots \otimes \mu_n$ on the space $X_1 \times
     \ldots \times X_n$. This means that we will be deriving a product form of the 
     invariance equation ($\ref{basicInvariance}$).
\end{enumerate}

After these considerations, the mathematics unfolds in much the same way
as before, with one exception. Invariant measures 
$\mu_1 \otimes \ldots \otimes \mu_n$ 
can exist only if a certain eigenvector condition involving the total measures
of the $\mu_i$ and the $\sigma^{}_{ij}$ is met (see Eq.~(\ref{eigCond}) below). 

Let $N$ and $(X_i,d_i)$, $i\in N$, be as above. For each 
$J=(j^{}_1,\dots,j^{}_n) \in N^n$, we define
$X^{}_J:= X_{j^{}_1} \times \ldots\times X_{j^{}_n}$ and adopt standard
multi-index notation, e.g.\ $x^{}_J=(x_{j^{}_1},\dots,x_{j^{}_n})$. 
In particular, $X^{}_N = X_1 \times \ldots  \times X_n$ in agreement
with our previous definition.  
We then define the metric $d^{}_J$ on $X^{}_J$ by $d^{}_J(x^{}_J,y^{}_J) = 
\sup_{k=1}^n d_{j^{}_k}(x_{j^{}_k},y_{j^{}_k})$.
For measures $\mu_i \in \mathcal{M}(X_i)$, $i\in N$, we write
$\mu^{}_J = \mu_{j^{}_1} \otimes \ldots  \otimes \mu_{j^{}_n} \in 
\mathcal{M}(X^{}_J)$ and 
${\rm d}\mu^{}_J= {\rm d}\mu_{j^{}_1} \dots {\rm d}\mu_{j^{}_n}$.

{}For each $(i,j) \in N\times N$, let $F_{ij}$  be a compact admissible family
of contractions $f\!:X_j \longrightarrow X_i$ (allowing, as above, the 
possibility that $F_{ij}$ is empty). We let $0 < r < 1$ be a uniform upper 
bound on the contractivity factors of all these mappings.

We let $F= \btimes_{i,j} \, F_{ij}$ be the product of all these spaces of maps, 
a typical element being a matrix of maps $\bs{f} = (f_{ij})$. For each such 
$\bs f$, and for all $J,K \in N^n$, let 
${f}^{}_{KJ}: X_J \longrightarrow X_K$ be given by
\begin{equation}
  {f}^{}_{KJ}(x_{j^{}_1}, \dots, x_{j^{}_n}) \; = \;
  (f_{k^{}_1 j^{}_1}(x_{j^{}_1}), \dots , f_{k^{}_n j^{}_n}(x_{j^{}_n})) \, .
\end{equation}
We write ${f}^{}_J$ for the special case 
${f}^{}_{NJ}\! : X^{}_J \longrightarrow X^{}_N$ and
$F^{}_J := \{{f}^{}_J \mid \bs{f} \in F\}$.
Note that now $f^{}_{KJ}.\mu^{}_J =
(f^{}_{k^{}_{1} j^{}_{1}}.\mu^{}_{j^{}_{1}})\otimes\ldots\otimes
(f^{}_{k^{}_{n} j^{}_{n}}.\mu^{}_{j^{}_{n}})$. Consequently,
$f^{}_{KJ}.\mu^{}_J \in \cM(X_K)$ and $f^{}_{J}.\mu^{}_J \in \cM(X_N)$.

We assume that each space $F_{ij}$ is equipped with a positive Borel measure
$\sigma^{}_{ij}$ and define $s^{}_{ij} := \sigma^{}_{ij}(F_{ij})$, or 
$s^{}_{ij}=0$ if $F_{ij}$ is empty. For each $J,K \in N^n$, we define the measure 
$\sigma^{}_J := \sigma^{}_{NJ} =
\sigma^{}_{1 j^{}_1}\otimes\ldots\otimes\sigma^{}_{n j^{}_n}$ 
and $s^{}_{KJ} := s^{}_{k_1j_1}\cdot \ldots \cdot s^{}_{k_nj_n}$.

The matrix $\bs{s} := (s_{ij})$ is a non-negative matrix. 
We now make the following compatibility assumption:
\begin{itemize}
\item[\bf CA] The total measures $m_i = \mu_i(X_i)$ of the $\mu_i$ are
   all (strictly) positive, and $\bs{m} := (m_1,\dots,m_n)^T$
   is an eigenvector of $\bs{s}$ with eigenvalue 1:
\begin{equation}\label{eigCond}
   \bs{s} \bs{m}\; = \;\bs{m}\, .
\end{equation}
\end{itemize}

{\sc Remark}: If $\bs{s}$ is non-negative, but $\bs{s} \bs{m} = \bs{m}$
for a vector $\bs{m}$ with all $m_i>0$ as we assume in {\bf CA}, the
eigenvalue 1 is also the spectral radius of $\bs{s}$ 
(see Appendix 2 of \cite{KT}, and Corollary 2.2 of it in particular) 
and thus its Perron-Frobenius (PF) eigenvalue.
Under the additional assumption of irreducibility of $\bs{s}$
(which we do not make!), $\bs{m}$ would be the unique PF eigenvector,
and primitivity of $\bs{s}$ would further imply that all other eigenvalues
of $\bs{s}$ were less than 1 in absolute value.

Let us also mention that there is no need to choose any particular
normalization here, but a convenient one would be 
$m^{}_N := m^{}_1\cdot \ldots \cdot m^{}_n = 1$.

Define $\cP^{\bs{m}}(X^{}_N)$ to be the space of all {\em product} measures
\index{measure!product} $\mu=\mu^{}_1 \otimes \ldots  \otimes \mu^{}_n$ 
where $\mu_i \in \cM_+^{m_i}(X_i)$, i.e.\ $\mu_i$ is a positive measure of 
total variation $\|\mu_i\| = \mu_i(X_i) = m_i$. 
For each $\bs{f} = (f_{ij}) \in F$, we define the operator
\begin{equation}
\begin{array}{rcl}
     \cA_{\bs{f}} \!: \; \cP^{\bs{m}}(X_N) & \longrightarrow & \cM(X_N)\\
             \mu\quad & \mapsto & 
      \cA_{\bs{f}}(\mu) \, := \, \sum_{J\in N^n} ( f^{}_J . \mu^{}_J )\, .
\end{array}
\end{equation}

{}For any  $\phi \in C(X_N,\RR)$, we have
\begin{equation}
   \cA_{\bs{f}}(\mu)(\phi) \; =  \;
   \sum_{J\in N^n}  \mu^{}_J (\phi\circ f^{}_J)  \, .
\end{equation}
In particular, if $\phi(x^{}_1,\dots, x^{}_n) = 
\phi^{}_1(x^{}_1)\cdot \ldots \cdot \phi^{}_n(x^{}_n)$
for some $\phi_i \in C(X_i,\RR)$, then this can be rewritten as
\begin{eqnarray*}
  \cA_{\bs{f}}(\mu)(\phi) & = &
\sum_{j^{}_1, \dots , j^{}_n} \mu_{j^{}_1}(\phi^{}_1 \circ f_{1 j^{}_1})
\cdot \ldots \cdot \mu_{j^{}_n}(\phi^{}_n \circ f_{n j^{}_n})\\
 & = &
   \Bigl(\,\sum_j (f_{1j}.\mu_j) (\phi^{}_1) \Bigr)
   \cdot \ldots \cdot 
   \Bigl(\,\sum_j (f_{nj}.\mu_j) (\phi^{}_n) \Bigr) \\
  & = &
   \sum_{J \in N^n}  \bigl(f_{1 j^{}_1}.\mu_{j^{}_1} \otimes \ldots \otimes 
           f_{n j^{}_n}.\mu_{j^{}_n} \bigr) (\phi) \, , 
\end{eqnarray*}
which, since the linear span of the product functions 
$\phi = (\phi^{}_1, \dots ,\phi^{}_n)$ is dense in $C(X_N,\RR)$, shows that 
\begin{equation}
   \cA_{\bs{f}}(\mu) \; = \; \sum_{J\in N^n} 
   f_{1 j^{}_1}.\mu_{j^{}_1} \otimes \ldots \otimes 
   f_{n j^{}_n}.\mu_{j^{}_n} \; = \;
   \sum_{J\in N^n} f^{}_{J} . \mu^{}_{J} \, .
\end{equation}

We define the {\em averaging } operator $\cA$ \index{averaging operator} 
on $\cP^{\bs{m}}(X^{}_N)$ by
\begin{eqnarray} \label{avOp}
 \cA(\mu) & := &  \int_F \cA_{\bs{f}}(\mu) \, {\rm d}\bs{\sigma}(\bs{f}) \nonumber\\
  & = & \int_F  \sum_{J\in N^n}  \Bigl(
        f_{1 j^{}_1}.\mu_{j^{}_1} \otimes \ldots  \otimes 
        f_{n j^{}_n}.\mu_{j^{}_n} \Bigr) \, 
        {\rm d}\sigma^{}_{1j^{}_1}(f_{1 j^{}_1}) \cdot\ldots\cdot  
        {\rm d}\sigma^{}_{nj^{}_n}(f_{n j^{}_n}) \nonumber \\
  & = & \Bigl(\,\sum_j \int_{F_{1j}}(f_{1j}.\mu_j)
        \, {\rm d}\sigma^{}_{1j}(f_{1j}) \Bigr) \cdot \ldots \cdot
        \Bigl(\,\sum_j \int_{F_{nj}}(f_{nj}.\mu_j)
        \, {\rm d}\sigma^{}_{nj}(f_{nj}) \Bigr)
        \nonumber \\
  & = & \sum_{J\in N^n} \int_{F^{}_J} 
        (f^{}_J.\mu^{}_J) \, {\rm d}\sigma^{}_J(f^{}_J) \, ,  
\end{eqnarray}
for $\mu = \mu^{}_1\otimes\ldots \otimes\mu_n$.
For any  $\phi \in C(X^{}_N,\RR)$, this reads:
\begin{equation}
    \int_F \cA_{\bs{f}}(\mu)(\phi)\, {\rm d}\bs{\sigma}(\bs{f}) \; = \; 
    \sum_{J\in N^n} \int_{F^{}_J} \mu_J (\phi\circ f_J) 
    \, {\rm d} \sigma^{}_J(f^{}_J) 
    \;\,\in\; \cM(X^{}_N) \, ,
\end{equation}
and, if $\phi(x_1,\dots, x_n) = \phi_1(x_1)\cdot \ldots \cdot \phi_n(x_n)$
for some $\phi_i \in C(X_i,\RR)$, this becomes
\begin{eqnarray*}
\lefteqn{\int_F \cA_{\bs{f}}(\mu)(\phi) \, {\rm d}\bs{\sigma}(\bs{f})}\\
& = &
  \Bigl(\,\sum_j \int_{F_{1j}}(f.\mu_j)
  \, {\rm d}\sigma^{}_{1j}(f) (\phi_1) \Bigr) \cdot \ldots \cdot
  \Bigl(\,\sum_j \int_{F_{nj}}(f.\mu_j)
  \, {\rm d}\sigma^{}_{nj}(f) (\phi_n) \Bigr) .
\end{eqnarray*}

This shows that the averaging process is of the sort envisaged in 
(\ref{basicInvariance}) and
\begin{equation}\label{productForm}
   \cA(\mu) \; = \;
   \Bigl(\, \sum_j \int_{F_{1j}} (f.\mu_j)
   \, {\rm d}\sigma^{}_{1j}(f) \Bigr) \otimes  \ldots \otimes 
   \Bigl(\, \sum_j \int_{F_{nj}} (f.\mu_j)
   \, {\rm d}\sigma^{}_{nj}(f) \Bigr) .
\end{equation}
Furthermore, consider
\begin{equation}\label{productForm-2}
   \sum_j \int_{F_{ij}} (f.\mu_j)
   \, {\rm d}\sigma^{}_{ij}(f)  \; \in \; \cM(X_i) .
\end{equation}
Since $f\in F_{ij}$ implies $\bs{1}^{}_{X_i} \circ f = \bs{1}^{}_{X_j}$,
(\ref{productForm-2}) satisfies
\begin{eqnarray*}
  \sum_j \int_{F_{ij}}(f.\mu_j) \, {\rm d}\sigma^{}_{ij}(f) 
  \, (\bs{1}^{}_{X_i}) & = & 
  \sum_j \int_{F_{ij}} \mu_j (\bs{1}^{}_{X_i} \circ f)
  \, {\rm d}\sigma^{}_{ij}(f) \\
  = \;\, \sum_j \, m^{}_j \int_{F_{ij}} {\rm d}\sigma^{}_{ij}(f)
  & = & \sum_j  s^{}_{ij} m^{}_j \;\, = \;\, m^{}_i \, .
\end{eqnarray*}
This shows  that our averaging operator stabilizes the space of product 
measures that we are considering:
\begin{prop}
   The averaging operator 
   $\cA := \int_F \cA_{\bs{f}} \, {\rm d}\bs{\sigma}(\bs{f})$ 
   of Eq.~{\rm (\ref{avOp})}
   maps the space $\cP^{\bs{m}}(X_N)$ of product measures with mass
   vector $\bs{m}$ into itself. \qed
\end{prop}

If $\phi \in C(X^{}_N,\RR)$ is a contraction, then so is 
$\phi \circ f^{}_J \!: X^{}_J \longrightarrow \RR$ for every $\bs{f} \in F$:
\begin{eqnarray*}
  \bigl|\phi \circ f^{}_J(x) - \phi \circ f^{}_J(x')\bigr|  & \le & 
  r^{}_\phi \, d^{}_N(f^{}_J(x), f^{}_J(x'))   \\
  & = & r^{}_\phi\mbox{ $\sup_{i\in N}$}\, 
        d_i\bigl(f_{ij_i^{}}(x_{j_i^{}}),f_{ij_i^{}}(x_{j_i^{}}')\bigr)\\
  & \le &  r^{}_\phi \, r \mbox{ $\sup_{i\in N}$}\, 
          d_{j_i^{}}(x_{j_i^{}} , x_{j_i^{}}') \\
  & = & r^{}_\phi \, r \, d^{}_J(x,x')
\end{eqnarray*}
for all $x,x' \in X^{}_J$.

Define a metric $L^{}_J$ on $\cM_+^{m^{}_J}(X^{}_J)$, the space of positive 
measures of total measure $m^{}_J := \prod_{i=1}^n m_{j_i^{}} > 0$, by
\begin{equation}
   L^{}_J(\mu,\nu) \; = \;
      \frac{1}{m^{}_J} \, \sup \, \bigl\{
      \,\left|\mu(\psi) - \nu(\psi)\right| \,\mid
      \psi\in {\rm Lip}(\le\! 1,X^{}_J,\RR) \bigr\}\, .
\end{equation}
This makes $\cM_{+}^{m^{}_J}(X^{}_J)$
into a complete metric space by Proposition \ref{metricL}.

Define $L$ on $\cP^{\bs{m}}(X^{}_N)$ by
\begin{equation}
   L(\mu, \nu) \; = \; \sup \, \{L^{}_J(\mu^{}_J,\nu^{}_J)
     \mid J\in N^n \}\, .
\end{equation}
\begin{prop}
  The operator $\cA \! : \cP^{\bs{m}}(X_N) \longrightarrow \cP^{\bs{m}}(X_N)$
  is a contraction with respect to the metric $L$,
  with contractivity factor at most $r$. 
\end{prop}
{\sc Proof:} Let $\mu,\nu \in \cP^{\bs{m}}(X_N)$.
In order to determine the $L^{}_K(\cA(\mu),\cA(\nu))$, we have to determine 
$\cA(\mu)^{}_K$ for any $K \in N^n$. Now, $\cA\mu = \cA(\mu)$ is a product 
measure, and from (\ref{productForm}) we find
\begin{eqnarray*}
(\cA\mu)^{}_K 
  & = & \Bigl(\,\sum_{j^{}_1} \int_{F^{}_{k^{}_1 j^{}_1}} \, 
         (f.\mu^{}_{j^{}_1})\, {\rm d}\sigma^{}_{k^{}_1 j^{}_1}(f) \Bigr) 
          \otimes  \ldots  \otimes 
        \Bigl(\,\sum_{j^{}_n} \int_{F^{}_{k^{}_n j^{}_n}} \,  
         (f.\mu^{}_{j^{}_n})\, {\rm d}\sigma^{}_{k^{}_n j^{}_n}(f) \Bigr) \\
  & = & \sum_J \int_{F^{}_{KJ}} (f.\mu^{}_J)\, {\rm d} \sigma^{}_{KJ}(f) \, .
\end{eqnarray*}
Thus
\begin{eqnarray*}
L^{}_K \bigl((\cA\mu)^{}_K, (\cA\nu)^{}_K \bigr)
  &=& L^{}_K\Bigl(
     \sum_J \int_{F^{}_{KJ}}(f.\mu^{}_J)\, {\rm d} \sigma^{}_{KJ}(f)\, , 
     \sum_J \int_{F^{}_{KJ}}(f.\nu^{}_J)\, {\rm d} \sigma^{}_{KJ}(f) \Bigr) \\
  &=& \frac{1}{m_K} \, \sup_{\psi}\,\left| \, \sum_J \int_{F^{}_{KJ}}
      \bigl((f.\mu_J) - (f.\nu_J) \bigr)(\psi) \,
      {\rm d} \sigma^{}_{KJ}(f) \, \right|\\
  &\le&  \frac{1}{m_K} \sum_J \, \sup_{\psi}  \int_{F^{}_{KJ}}
      \bigl|\,(\mu^{}_J - \nu^{}_J) (\psi \circ f) \bigr| 
      \, {\rm d} \sigma^{}_{KJ}(f)  \\
  &\le&  \frac{r}{m_K} \sum_J \, \int_{F^{}_{KJ}} \sup_{\psi}\, \left|\, 
         (\mu^{}_J - \nu^{}_J) (r^{-1}\psi \circ f)\right| 
         \, {\rm d} \sigma^{}_{KJ}(f)  \\
  &\le&  \frac{r}{m_K} \sum_J \, 
         \int_{F^{}_{KJ}} m^{}_J L^{}_J (\mu^{}_J,\nu^{}_J)\, 
         {\rm d} \sigma^{}_{KJ}(f)  \\
  &\le&  \frac{r L(\mu,\nu)}{m_K} \sum_J \, m^{}_J \!\int^{}_{F^{}_{KJ}}  
         \!\! {\rm d} \sigma^{}_{KJ}(f) \; = \;
         \frac{r L(\mu,\nu)}{m_K} \sum_J \, s^{}_{KJ} m^{}_J \, ,
\end{eqnarray*}
where $\psi$ runs through ${\rm Lip}(\le\! 1, X^{}_K,\RR)$. {}Finally, 
\begin{eqnarray*}
    L(\cA \mu, \cA \nu ) 
  & = & \sup_K L^{}_K \bigl((\cA \mu)^{}_K, (\cA \nu)^{}_K\bigr) 
        \;\; \le \;\; r L(\mu,\nu) 
        \sup_K \frac{1}{m_K} \sum_J s^{}_{KJ}m^{}_J  \\
  & = &  r L(\mu,\nu) \sup_K \frac{1}{m^{}_K}\prod_{i=1}^{n}
         \sum_{j_i^{}=1}^n s^{}_{k_i^{} j_i^{}} m^{}_{j_i^{}}    \\
  & = &  r L(\mu,\nu) \sup_K \frac{1}{m^{}_K} \prod_{i=1}^n m^{}_{k_i^{}} 
         \; = \;\,  r L(\mu,\nu) \, ,
\end{eqnarray*}
where we have used (\ref{eigCond}) in the last line.
We have thus established the following result.
\begin{theorem}\label{com-met-spa}
  Let $X^{}_1, \dots, X^{}_n$ be compact metric spaces and, for each
  pair $(i,j)\in N\times N$, 
  let $F_{ij}$ be a compact admissible family of contractions
  $($possibly empty$)$ from $X_j$ to $X_i$. Assume that each $F_{ij}$
  is equipped with a positive Borel measure $\sigma_{ij}$
  and define $s_{ij} = \sigma_{ij}(F_{ij})$, with
  $s_{ij}:= 0$ if $F_{ij}$ is empty.  Assume that
  $\bs{s}:= (s_{ij})$ has a positive $1$-eigenvector
  $\bs{m} = (m_1, \dots, m_n)^T$. 

  Then there exists a unique family of measures 
  $\omega_i \in \cM(X_i)^{m_i}_+$ which satisfy
\begin{equation} \label{multi-sefl-sim-meas}
  \omega_i \; = \; \sum_{j=1}^n \int_{F_{ij}} (f.\omega_j)
  \, {\rm d} \sigma^{}_{ij}(f) \; , \quad i\in N \, .  
 \end{equation}
The support of $\omega$ is contained in $W$, the attractor of
{\em (\ref{multiAttractor})}.  \qed
\end{theorem}

We call $\omega = \omega^{}_1 \otimes \ldots \otimes \omega^{}_n$ the 
$(F,\bs{\sigma},\bs{m})$-invariant measure on $X_N$, or simply
$(F,\bs{\sigma})$-invariant measure, if $\bs{m}$ is 
understood from the context.

\section{Multi-component families of affine mappings}\label{mul-comp}

In this Section, $H$ is an LCAG and $A$ is an automorphism of $H$ with
modulus $\alpha > 1$. The Haar measure on $H$ is denoted by $\theta$.
$H$ is assumed to be complete with respect to a metric $d$, relative to
which $A$ is a contraction.

We assume that we are given $n$ copies of $H$, which we call $H_1,\dots,H_n$,
and compact families $F_{ij}$, $1\le i,j\le n$, of affine mappings
$f_{ij}\!: H_j\to H_i$, all of the form $f_{ij}(x)=Ax+u^{}_{ij}$ with 
$u^{}_{ij}\in H_i$.
Our objective is to understand the multi-component system formed by
$$ H^n \; = \; H_1\times\ldots \times H_n \; = \; H\times\ldots \times H $$
and the admissible family of contractions $F = \btimes F_{ij}$ in a way
that parallels our previous analysis in Section \ref{affine-maps}.

By Proposition \ref{lip-contra}, $F$ has a unique attractor 
$W=W_1\times\ldots\times W_n\in (\cK H^n)$. We wish to describe the unique
$F$-self-similar measure $\omega=\omega^{}_1\otimes\ldots\otimes\omega^{}_n$
on $H^n$ with respect to a system $\bs{\sigma}=(\sigma^{}_{ij})$ of measures
on $F$.

The mapping $F_{ij}\to H_i$ defined by $f_{ij}\mapsto f_{ij}(0) = u_{ij}$ 
is continuous
and produces a homeomorphism between $F_{ij}$ and a compact subset 
$F_{ij}':= F_{ij}(0)$ of $H_i$.

We assume that each compact space $F_{ij}$ is supplied with a positive
regular Borel measure $\sigma^{}_{ij}$, supported on $F_{ij}$ with
$s^{}_{ij}:=\sigma^{}_{ij}(F_{ij})$. We may identify $\sigma_{ij}$
with a regular Borel measure on $H_i$ supported on $F_{ij}'$.
It is understood that $F_{ij}$ may be empty, in which case $s^{}_{ij}:=0$.
Furthermore, we assume the existence of a mass vector
$\bs{m}=(m^{}_1,\dots,m^{}_n)^T > 0$
satisfying the compatibility assumption {\bf CA}, i.e.\ $\bs{s}\bs{m}=\bs{m}$.

We define $X_i\in\cK H_i$, $i=1,\dots,n$, to be compact subspaces with the
following properties:
\begin{enumerate}
\item $X_1\times\ldots \times X_n$ is invariant under the family of mappings
      $F$;
\item $\sum_{k=0}^{\ell} A^k \,\textrm{supp}(\bs{\sigma}) \subset
      X_1\times\ldots \times X_n$, for all $\ell\ge 0$.
\end{enumerate}
It is easy to see that such sets exist because the mappings of $F$ and the 
automorphism $A$ are all contractive. The $F$-invariance already forces
$W_i\subset X_i$.

Let notation be as in Section \ref{mcms}, so $X^{}_N=X_1\times\ldots \times X_n$ 
and $\cP^{\bs{m}}(X^{}_N)$ is the space of all product measures 
$\mu^{}_1\otimes\ldots \otimes\mu^{}_n$ on $X^{}_N$ for which 
$\mu^{}_i\in\cM_{+}^{m_i}(X_i)$. We know that the averaging operator $\cA$
of (\ref{avOp}) is a contraction on $\cP^{\bs{m}}(X^{}_N)$.

Let $E_i\subset X_i$, $i\in N$, be measurable sets. 
Then, from (\ref{productForm}), we obtain
\begin{eqnarray}
  \cA\mu (E_1\times\ldots \times E_n)
  & = & \prod_{i=1}^{n} \sum_{j} \int_{F_{ij}}
  (f.\mu^{}_{j}) (E_i)\, {\rm d}\sigma^{}_{ij} (f)\nonumber \\
  & = & \prod_{i=1}^{n} \sum_{j} \int_{F_{ij}}
  \mu^{}_{j} (f^{-1} (E_i))\, {\rm d}\sigma^{}_{ij} (f) \nonumber \\
  & = & \prod_{i=1}^{n} \sum_{j} \int_{F_{ij}'}
  (A.\mu^{}_{j})(E_i - u) 
  \, {\rm d}\sigma^{}_{ij}(u) \nonumber \\
  & = & \prod_{i=1}^{n} \sum_{j} (\sigma^{}_{ij}*A.\mu^{}_{j})(E_i)\, .
\end{eqnarray}
Thus
$$ \cA\mu \; = \; \Bigl(\,\sum_j \sigma^{}_{1j} * A.\mu^{}_j \Bigr)
   \otimes\ldots\otimes \Bigl(\,\sum_j \sigma^{}_{nj} * A.\mu^{}_j \Bigr).$$
Adopting matrix notation, with $\mu=\mu^{}_1\otimes\ldots \otimes\mu^{}_n$
written as $(\mu^{}_1,\dots,\mu^{}_n)^T$ and $\bs{\sigma}=(\sigma^{}_{ij})$,
this reads
\begin{equation}\label{matrix-form}
   \cA\mu \; = \; \bs{\sigma} * A.\mu
\end{equation}
where $A.\mu := (A.\mu^{}_1,\dots,A.\mu^{}_m)^T$. If we now define
$A.\bs{\sigma}:=(A.\sigma^{}_{ij})$, we can iterate (\ref{matrix-form}). 
Observing Proposition \ref{facts}(3), we obtain
\begin{equation}\label{conv-iter}
   \cA^{\ell} \mu \; = \; \bs{\sigma} * A.\bs{\sigma} * \ldots * 
   A^{\ell -1}.\bs{\sigma} * A^{\ell}.\mu \, .
\end{equation}

We now proceed as in Section \ref{affine-maps} to define a suitable sequence 
of measures $\bigl(\omega^{(\ell)}\in\cP^{\bs m}(X)\bigr)_{\ell\ge 0}$.  
First, let
\begin{equation}\label{pointmeasure}
  \omega^{(0)} \; := \; \delta^{\bs{m}} \; = \; 
          (m^{}_1\delta^{}_0, \dots, m^{}_n \delta^{}_0)^T 
\end{equation}
where $\delta_0$ is the unit point measure supported at $\{0\}$.
Clearly, $\delta^{\bs{m}}\in\cP^{\bs m}(X)$, but since
$A$ is a contraction, we also have $A.\delta^{\bs{m}}=\delta^{\bs{m}}$, and 
$A^{\ell}.\mu\to\delta^{\bs{m}}$ as $\ell\to\infty$, for any
$\mu\in\cP^{\bs m}(X)$. Define iteratively, as before, 
$\omega^{(\ell+1)} = \cA\omega^{(\ell)}$ for $\ell\ge 0$. Then 
$$ \omega^{(\ell+1)} \; = \; \bs{\sigma} * A.\bs{\sigma} * \ldots *
   A^{\ell}.\bs{\sigma} * \delta^{\bs{m}} \, . $$
We have $\textrm{supp}(\omega^{(\ell)})\subset\sum_{k=0}^{\ell}
A^k\,\textrm{supp}(\bs{\sigma}* \delta^{\bs{m}} )$ 
and $\cA\cP^{\bs{m}}(X)\subset\cP^{\bs{m}}(X)$, so we know,
since $\delta^{\bs{m}}\in\cP^{\bs{m}}(X)$, that
$\omega^{(\ell)}\in\cP^{\bs{m}}(X)$. Consequently, $\omega^{(\ell)}$
vaguely converges, as $\ell\to\infty$, to the unique 
$(F,\bs{\sigma},\bs{m})$-self-similar measure
$\omega =\omega^{}_1\otimes\ldots \otimes\omega^{}_n\in\cP^{\bs{m}}(X)$.
Since $\textrm{supp}(\omega) \subset W$ by Theorem \ref{com-met-spa},
we find that $\omega\in\cP^{\bs{m}}(W)$. To summarize:

\begin{prop}\label{aff-inv-meas}
  Let $H$ be an LCAG which is a complete metric space. Let $A$ be a contractive
  automorphism on $H$ and let $F_{ij}$, $1\leq i,j\leq n$, with attractor
  $W$, be compact admissible families of affine maps on $H$, all of the form
  $x\mapsto Ax+v$ with $v\in H$. Let $\sigma^{}_{ij}$ be a positive regular Borel
  measure on $F_{ij}$, identified with a regular Borel measure on $H$ supported
  on $F_{ij}(0)$ $($with $\sigma^{}_{ij}:=0$ if $F_{ij}=\emptyset \,)$. 
  Let $\bs{s}=(s^{}_{ij})=(\sigma^{}_{ij}(H))$ and, finally, let
  $\bs{m}=(m^{}_1,\dots ,m^{}_n)^T >0$ satisfy $\bs{s}\bs{m}=\bs{m}$. Then
\begin{enumerate}
\item $\vphantom{\Bigl(\conv_1^2\Bigr)}
      \omega \, = \, \bigl( \conv_{\ell=0}^{\infty}\, A^{\ell}.\bs{\sigma}
      \bigr) * \delta^{\bs{m}}$ is the unique 
      $(F,\bs{\sigma},\bs{m})$-self-similar measure of 
      {\em (\ref{multi-sefl-sim-meas})}, with $\omega\in\cP^{\bs{m}}(H^n)$,
      ${\rm supp}(\omega)\subset W$, and $\delta^{\bs{m}}$ as
      in {\rm (\ref{pointmeasure})}.
\item $\vphantom{\Bigl(\conv_1^2\Bigr)}
      \widehat{\omega} \, = \, \left(\,\prod_{\ell=0}^{\infty} \,
      (A^T)^{-\ell}.\widehat{\bs{\sigma}}\,\right) \,\bs{1}^{\bs{m}}$, where
      $\bs{1}^{\bs{m}}=(m^{}_1 \bs{1}^{}_H,\dots,m^{}_n \bs{1}^{}_H)^T$,
      the convergence of the product being uniform on compact sets.
      \qed
\end{enumerate}
\end{prop}

If we assume that the measures $\sigma^{}_{ij}$ are absolutely continuous
with respect to Haar measure $\theta$ on $H$, then 
$\sigma^{}_{ij}=h^{}_{ij} \theta$ where $h^{}_{ij}\in L^1(H)$ due to the
Radon-Nikodym theorem. In particular, we have 
$\textrm{supp}(h^{}_{ij})\subset F_{ij}(0)\subset H$,
$h^{}_{ij} \ge 0$ and $\|h^{}_{ij}\|^{}_1 = 
\int_H h^{}_{ij}\, {\rm d}\theta =s^{}_{ij}$, for all $1\leq i,j\leq n$. Then
\begin{eqnarray*}
  \omega^{(\ell+1)} & = & \bs{\sigma} * A.\bs{\sigma} * \ldots *
     A^{\ell}.\bs{\sigma} * \delta^{\bs{m}} \\
  & = & \bs{h}\Theta * A.(\bs{h}\Theta) * \ldots * A^{\ell}.(\bs{h}\Theta) 
     \, * \, \delta^{\bs{m}} \\
  & = & \bs{h}\Theta * \alpha(A.\bs{h})\Theta * \ldots *
     \alpha^{\ell}(A^{\ell}.\bs{h})\Theta \, * \, \delta^{\bs{m}} \\
  & = & \bs{h} * \alpha(A.\bs{h})*\ldots *\alpha^{\ell}(A^{\ell}.\bs{h}) \,
     (\Theta * \delta^{\bs{m}})
\end{eqnarray*}
where $\bs{h}=(h^{}_{ij})$ and $\Theta = \textrm{diag}(\theta,\dots,\theta)$
is a diagonal matrix. Thus we have $\Theta * \delta^{\bs{m}} =
(m^{}_1\theta,\dots,m^{}_n\theta)^T$ and
\begin{equation}\label{star}
   \omega^{(\ell+1)} \; = \; \left(
   \overset{\ell}{\underset{k=0}{\conv}} \,
   \alpha^{k} (A^{k}.\bs{h}) \right)
   \,(m^{}_1\theta,\dots,m^{}_n\theta)^T 
   \;\, \in \;\, \cP^{\bs{m}}(H^n) \, .
\end{equation}
Vague convergence of this sequence is clear, but the results of Section
\ref{mcms} suggest that we can expect more. However, $\|.\|$-convergence is 
technically
more involved here. Let us thus first postpone this question and state
first the result on the self-similar functions.

\begin{prop}\label{cont-inv-meas}
Let notation and assumptions be as in Proposition {\rm\ref{aff-inv-meas}}, and suppose
that the measures $\sigma^{}_{ij}=h^{}_{ij}\theta$ are absolutely continuous
with respect to Haar measure $\theta$. Assume that the convolution in
{\rm (\ref{star})}, as $\ell\to\infty$, converges also in the $\|.\|$-topology.
Then there is a unique vector $\bs{g}=(g^{}_1,\dots,g^{}_n)^T$ of 
non-negative functions in $L^1(H)$ that satisfies
\begin{enumerate}
\item $\bs{g}\, = \, \bigl(\conv_{\ell=0}^{\infty} \,
      \alpha^{\ell} (A^{\ell}.\bs{h})\bigr) \, \bs{m}$
\item $\vphantom{\Bigl(\conv_{\ell=0}^{\infty}\Bigr)}
      g^{}_i (x) \, = \, \sum_{j=1}^{n} \int_H h^{}_{ij} (x-v) \,
      g^{}_{j}(A^{-1} v) \, {\rm d}\theta (v)\,$, $i=1,\dots,n$, with
      normalization $\int_H g^{}_i\, {\rm d}\theta = m^{}_i$ and
      $\,{\rm supp}(g^{}_i)\subset W_i$.
\end{enumerate}
Furthermore, if the $h^{}_{ij}$ are functions in $L^1(H)\cap L^{\infty}(H)$,
then the functions $g^{}_i$ are continuous on $H$.
\end{prop}
{\sc Proof}:
Part 1 is a direct reformulation of Eq.~(\ref{star}), and Part 2 is a
component-wise recounting of Part 1. The continuity follows from the
properties of the convolution product, as in Proposition \ref{single-conv}(4). \qed   

\smallskip
Infinite convolution products like that of Proposition \ref{cont-inv-meas}(1) 
also appear in the context of matrix continuous refinement operators.
These are introduced in \cite{JL} (with $H$  being $\RR^n$).
In our paper \cite{BMmulti}, we relied on the results of \cite{JL} for the
existence of our self-similar densities. However, the methods of
\cite{JL} are from functional analysis and do not lend themselves to the 
general measure theoretic situation that we are trying to address here.

Let us come back to the convergence issue in Eq.~(\ref{star}). Unlike the
situation in Section \ref{self-sim-functions}, with Lemma \ref{approx-unit} 
and Proposition \ref{norm-conv}, the $\|.\|$-convergence in (\ref{star}) is 
not entirely automatic. Note that the vector notation for the measures is
handy for the formulation of the iteration, but it still
represents a product measure. We are interested in the $\|.\|$-convergence
of the sequence of product measures (\ref{star}).
For this, it is sufficient, but not necessary, that the sequence
$K_{\ell} = \conv_{k=0}^{\ell} \, \alpha^k (A^k . \bs{h})$, seen as a sequence 
of linear operators, converges in the operator norm. However, for fixed $i,j$,
$\| (K_{\ell})_{ij} \|^{}_1 = (\bs{s}^{\ell+1})_{ij}$, and convergence of this,
for $\ell\to\infty$, does not follow from our general assumptions on the
matrix $\bs{s}$, see the Remark after Eq.~(\ref{eigCond}), because we did not
assume primitivity of $\bs{s}$.

Nevertheless, there is an analogue of Proposition \ref{norm-conv} which we will
now derive. To this end, define $\bs{h}^{(k)} = \alpha^k (A^k . \bs{h})$ for
$k\ge 0$. In particular, $\bs{h}^{(0)}=\bs{h}=(h_{ij})$, 
which is a matrix of functions
in $L^1(H)$, and also each $(\bs{h}^{(k)})_{ij}$ is a non-negative $L^1$-function
of norm $s^{}_{ij}$. Recall that (\ref{star}) means $\omega^{(\ell)} =
f^{(\ell)}_1 \theta \otimes \ldots \otimes f^{(\ell)}_n \theta$,
with $L^1$-functions $f^{(\ell)}_i = \sum_{j}
(\conv_{k=0}^{\ell-1}\,\bs{h}^{(k)})^{}_{ij} \, m^{}_j$ of norm
$\| f^{(\ell)}_i \|^{}_1 = \sum_j (\bs{s}^{\ell})^{}_{ij}\, m^{}_{j} = m^{}_{i}$.
Consequently, showing that (\ref{star}) converges also in the $\|.\|$-topology
means showing that $f^{(\ell)}_i$ converges in $L^1(H)$ for each $i$ as 
$\ell\to\infty$.

Fix $\varepsilon>0$, and let $V_{ij}$ be the corresponding neighbourhood for 
the function $h_{ij}$ according to Lemma \ref{approx-unit}. Let
$V=\bigcap_{i,j} V_{i,j}$ and choose an integer $M$ such that, for all
$k\ge 0$ and all $i,j$, the non-negative $L^1$-function
$(\conv_{\ell=M}^{M+k}\,\bs{h}^{(\ell)})_{ij}$, of norm
$(\bs{s}^{k+1})_{ij}$, has support inside $V$. Such an $M$ clearly exists.
With Lemma \ref{approx-unit}, we then find, for all $i',j'$ simultaneously,
the approximation formula
$$\bigl\| (\conv_{\ell=M}^{M+k}\,\bs{h}^{(\ell)})_{ij} * 
h_{i'j'} - (\bs{s}^{k+1})_{ij} \, h_{i'j'} \bigr\|_{1} 
\; \le \; (\bs{s}^{k+1})_{ij} \, \varepsilon \, .$$
Note that this formulation remains valid even in the limiting case that
$(\bs{s}^{k+1})_{ij}$ happens to vanish.

Let now $n,m\ge M$ and define 
$u^{}_{ij} = (\conv_{\ell=M}^{n}\,\bs{h}^{(\ell)})^{}_{ij}$ and
$v^{}_{ij} = (\conv_{\ell=M}^{m}\,\bs{h}^{(\ell)})^{}_{ij}$. 
Then, we can calculate as follows
\begin{equation*}
\begin{split}
  \bigl\| f^{(n)}_{i} - f^{(m)}_{i} \bigr\|_1
  & \; = \; \bigl\| \, \mbox{$\sum_{k,\ell,j}$} \,
   (\bs{h}^{(1)} * \ldots * \bs{h}^{(M-1)})^{}_{k\ell}
   * h^{}_{ik} * (u^{}_{\ell j} - v^{}_{\ell j} )\, m^{}_{j} \bigr\|_1 \\
  & \; \le \; \mbox{$\sum_{k,\ell}$} \, \left( 
   \bigl\| (\bs{h}^{(1)} * \ldots * \bs{h}^{(M-1)})^{}_{k\ell} \bigr\|_1
   \cdot \bigl\| \, \mbox{$\sum_j$} \, h^{}_{ik}
   * (u^{}_{\ell j} - v^{}_{\ell j} )\, m^{}_{j} \bigr\|_1 \right) \\
  & \; \le \; \mbox{$\sum_{k,\ell}$}\, (\bs{s}^{M-1})^{}_{k \ell} \;
   \bigl\| \, \mbox{$\sum_j$}\,  h^{}_{ik} * (u^{}_{\ell j} - v^{}_{\ell j})
   m^{}_j \bigr\|_1
\end{split}
\end{equation*} 
where we have used that convolution on the level of functions is commutative.
Observe next that
\begin{eqnarray*}
\lefteqn{\bigl\| \,\mbox{$\sum_j$}\, h^{}_{ik} * (u^{}_{\ell j} - v^{}_{\ell j}) 
    \, m^{}_{j} \bigr\|_1 \; \le \;  \bigl\|\,h^{}_{ik} \,\mbox{$\sum_j$}\,
    \bigl( (\bs{s}^{n-M+1})^{}_{\ell j} - 
    (\bs{s}^{m-M+1})^{}_{\ell j} \bigr) m^{}_{j} \bigr\|_1} \\
 & & + \; \mbox{$\sum_j$}\, \Bigl(
    \bigl\| h^{}_{ik} * u^{}_{\ell j} - h^{}_{ik} \,
    (\bs{s}^{n-M+1})^{}_{\ell j} \bigr\|_1 +
    \bigl\| h^{}_{ik} * v^{}_{\ell j} - h^{}_{ik} \,
    (\bs{s}^{m-M+1})^{}_{\ell j} \bigr\|_1 \Bigr)\, m^{}_{j} \\
 & \le & 
 \varepsilon \bigl(\, \mbox{$\sum_j$}\,
   (\bs{s}^{n-M+1})^{}_{\ell j} \, m^{}_{j} + \mbox{$\sum_j$}\,
   (\bs{s}^{m-M+1})^{}_{\ell j} \, m^{}_{j} \bigr)
   \;\, = \;\, 2 m^{}_{\ell} \, \varepsilon
\end{eqnarray*}
where we have used the above approximation formula and the equation
$\bs{s}\bs{m}=\bs{m}$. This finally gives
$$
 \bigl\| f^{(n)}_{i} - f^{(m)}_{i} \bigr\|_1
    \; \le \; 2\,\varepsilon\, \mbox{$\sum_{k,\ell}$}\, (\bs{s}^{M-1})^{}_{k \ell}
   \, m^{}_{\ell} 
   \; = \; 2 \,(m^{}_1 + \ldots + m^{}_n) \,\varepsilon \, , $$
independently of $i$. This shows that all sequences $(f^{(n)}_{i})_{n\ge 0}$
are Cauchy, and we have thus established the expected analogue of
Proposition \ref{norm-conv} and strengthening of Proposition \ref{cont-inv-meas}:
\begin{prop} 
Let notation and assumptions be as in Proposition {\rm\ref{aff-inv-meas}}, and 
suppose that the measures $\sigma^{}_{ij}=h^{}_{ij}\theta$ are absolutely continuous
with respect to Haar measure $\theta$. Then, the sequence of product measures in
{\rm (\ref{star})}, as $\ell\to\infty$, converges not only vaguely, but
also in the $\|.\|$-topology of $\cP^{\bs{m}}(H\times\ldots\times H)$.
\qed
\end{prop}

\section{Model sets and Weyl's Theorem}\label{model-weyl}

To link our previous analysis to quasicrystals, let us now summarize some
of the key ingredients to their mathematical description.
A {\em cut and project scheme} \index{cut and project scheme} consists of
the following set of data:
\begin{itemize}
\item a real space $\RR^m$
\item a locally compact Abelian group $H$
\item a lattice $\tilde{L}\subset\RR^m\times H$
\end{itemize}
which satisfies the following properties. If $\pi$ and $\pi^{}_H$ are the
natural projections of $\RR^m\times H$ onto $\RR^m$ and $H$, respectively,
then
\begin{itemize}
\item $\pi|^{}_{\tilde{L}}$ is one-to-one.
\item $\pi^{}_H(\tilde{L})$ is dense in $H$.
\end{itemize}
This is summarized in the following diagram.
\begin{equation}\label{diagram}
\begin{array}{ccccc}
   \RR^m & \stackrel{\pi}{\longleftarrow} &
   \RR^m\times H & \stackrel{\pi^{}_H}{\longrightarrow} & H \\
%   &&&&\\ && \uparrow && \\ &&&&\\
   && \mbox{\raisebox{-1ex}{$\cup$}} && \\
     \mbox{\raisebox{1.5ex}{\footnotesize 1-1}}
     \hspace*{-3ex}& 
     \mbox{\raisebox{2ex}{$\nwarrow$}}   & 
     \mbox{\raisebox{-1ex}{$\tilde{L}$}} &
     \mbox{\raisebox{2ex}{$\nearrow$}}   & \hspace*{-2ex}
     \mbox{\raisebox{1.5ex}{\footnotesize dense}}
\end{array}
\end{equation}
To say that $\tilde{L}$ is a {\em lattice} \index{lattice} in $\RR^m\times H$ 
means that $\tilde{L}$ is a discrete subgroup of $\RR^m\times H$ such that
$(\RR^m\times H)/\tilde{L}$ is compact.

We set $L=\pi(\tilde{L})$, a subgroup of $\RR^m$, and define
the {\em star map} \index{star map} $(.)^*\!: L\to H$ by
$x^* = \pi^{}_H\circ\bigl(\pi|^{}_{\tilde{L}}\bigr)^{-1} (x)$.
Although $(.)^*$ is a group homomorphism, it has, in general, no natural
extension to $\RR^m$ and, indeed, it is typically totally discontinuous in
the topology on $L$ induced by $\RR^m$. In fact, it is this property
that makes it useful!

Given any subset $U\subset H$, we define
\begin{equation}\label{premodelset}
   \Lambda(U) \; := \; \{ x\in L\mid x^*\in U\} \; = \;
   \{\pi(\tilde{x})\mid\tilde{x}\in\tilde{L}\, ,\;
   \pi^{}_H(\tilde{x})\in U\} \; \subset \; L \; \subset \; \RR^m\, .
\end{equation}
A set $\Lambda\subset\RR^m$ is a {\em model set} \index{model set}
relative to (\ref{diagram}), if $\Lambda=\Lambda(W)$ for some 
$W\subset H$ that is compact and equals the closure of its 
non-empty interior\footnote{There are variations on the exact conditions 
imposed on $W$ depending on the delicateness of the results required. 
Our assumptions imply the Delone property of $\Lambda$, and
are rather convenient for many other purposes. 
There is still something unsatisfying about
our present understanding of model sets. To say that $\Lambda\subset\RR^m$
is a model set is to say that it arises from some cut and project scheme.
But we still lack a useful direct characterization of such sets, compare
\cite{Martin}.}.

Model sets have remarkable properties that make them important objects of
study in the theory of mathematical quasicrystals. We refer the reader
to \cite{MB,Moody,P} and references therein for more details, but we mention
here a few of the key points:
\begin{enumerate}
\item If $\Lambda\subset\RR^m$ is a model set, then it is a {\em Delone set}
      \index{Delone set} in $\RR^m$, i.e.\ $\Lambda$ is both uniformly 
      discrete and relatively dense. 
\item Generically, model sets have no translational symmetries, although
      they certainly still have a high degree of long-range order.
\item If $W$ is {\em Riemann measurable}, \index{Riemann measurable}
      i.e.\ if $\partial W$ has vanishing Haar measure in $H$, then
      $\Lambda$ has a well-defined density (see Proposition 
      \ref{SchlottmannThm} below).
\item If $W$ is {\em Riemann measurable}, then $\Lambda$ is
      pure point diffractive \cite{Hof-diff,Martin2}.
\end{enumerate}

Model sets appeared very early in the theory of quasicrystals (under the
name of cut and project sets, see \cite{MB} and references given there), but 
originally only with the internal group
$H$ being another real space. However, model sets had been defined much
earlier and in full generality in a totally different context by Y.~Meyer
\cite{Meyer}. Recent papers \cite{Moody,BMS,Martin2} show that the more general
setting is completely relevant in the mathematical theory of quasicrystals,
aperiodic tilings and substitution systems, both geometric and algebraic.

A key feature of a model set is that it puts together a discrete geometric
object, $\Lambda(W)$, with a relatively compact set $W\subset H$ on which we 
can use the powerful array of tools from analysis on locally compact Abelian 
groups. The essential mathematical link is Weyl's Theorem \index{Weyl's Theorem}
on uniform distribution \cite{Weyl} which connects densities on $\Lambda(W)$
to measures on $W$. We refer to this theory in the more general setting of
LCAG's, see \cite[Ch.\ 4.4]{KN} for background material.

In its usual form, Weyl's Theorem is stated for real spaces, 
but it works at the level of locally compact Abelian groups, too. Here we
establish the theorem in this more general setting. The basis of the theorem 
in the context of model sets is the phenomenon of `uniformity of projection', 
which we state here in the generality that we will need. In fact, our proof 
demonstrates that Weyl's Theorem, in this context, is actually equivalent
to the uniformity of projection. For an approach to uniformity of projection 
via ergodic theory, see \cite{Hof,Martin2}.

In the following two Propositions, it is understood that a cut and project 
scheme according to (\ref{diagram}) has been given. In addition, let
$\tilde{\theta}$ denote the product measure on $\RR^m\times H$, formed from 
Lebesgue measure on $\RR^m$ and our fixed Haar measure $\theta$ on $H$.
% Let us use the symbol $\tilde{\theta}$ also for the induced measure on
% $T:=(\RR^m\times H)/\tilde{L}$ which can be interpreted as the
% fundamental domain of the lattice $\tilde{L}$.
% We use $|\tilde{L}|:=\tilde{\theta}(T)$ for its volume in the sequel. 
Let now $T$ be any measurable fundamental domain of $\RR^m\times H$ with
respect to the action of its discrete subgroup $\tilde{L}$, and define 
$|\tilde{L}|:=\tilde{\theta}(T)$ as its volume. Note that the value of 
$|\tilde{L}|$ does not depend on the actual choice of $T$. Its meaning
really is the averaged number of lattice points per unit volume
(in Haar measure).

\begin{prop} \label{SchlottmannThm}
{\rm (Schlottmann \cite{Martin})} 
   Let a cut and project setup according to diagram {\rm (\ref{diagram})}
   be given, with $|\tilde{L}|$ as described above.
   Let $U \subset H$ be totally bounded and Riemann measurable 
   {\rm (}i.e.\ $U$ measurable with $\theta(\partial U)=0${\rm )}. Then 
$$ \lim_{r \rightarrow \infty} \frac{1}{{\rm vol}\, B_{r}(0)} 
   \left( \sum_{x \in \Lambda(U) \cap B_{r}(a)}1 \right) \; = \;
   \frac{\theta(U)}{|\tilde{L}|} \; . $$
Furthermore, the limit is uniform in $a$. \qed
\end{prop}
This limit is called the {\em density}, $\textrm {den}\, (\Lambda(U))$, of 
$\Lambda(U)$. \index{density} Note that the totally bounded set $U$ in this 
Proposition need not be closed. It is only demanded that its boundary has
vanishing Haar measure. If $U$ itself is of measure 0, the density of
$\Lambda(U)$ vanishes.

{\sc Remark}: There is another way to explain the meaning of $|\tilde{L}|$
which is perhaps more natural from the group theoretic point of view.
Consider the factor group $T':=(\RR^m\times H)/\tilde{L}$ (which is compact)
and let $\mu$ be its normalized Haar measure. If now $f$ is a continuous
function on $\RR^m\times H$ with compact support, define a new function by
$F(x) = \sum_{u\in\tilde{L}} f(x+u)$. So, $F$ results from $f$ by averaging
over the canonical Haar measure of the lattice $\tilde{L}$, which is counting 
measure. The function $F$ can then also be viewed as a
function on $T'$, and we can determine its integral, $\mu(F)$. If we then
define a new measure on $\RR^m\times H$ by $\tilde{\theta}'(f) := \mu(F)$,
it is another Haar measure on $\RR^m\times H$, and we must have
$\tilde{\theta}' = c \, \tilde{\theta}$. The constant $c$ is nothing
but $|\tilde{L}|$, see \cite[Ch.\ XIV.4]{Dieu} for background material.

\begin{theorem}  {\rm (Weyl's Theorem for general model sets)}
   \label{weyl-thm} \index{Weyl's Theorem}
   Let $\Lambda=\Lambda(W)$ be a model set in the above sense,
   with compact, Riemann measurable $W\subset H$.
   Let $f\!:H \rightarrow \mathbb R$ be continuous with 
   ${\rm supp}\,(f)\subset W$.
   Let $p\!:L \rightarrow \RR$ be defined by $p(x) = f(x^*)$.   
   Then 
$$ \lim_{r \rightarrow \infty} \frac{1}{{\rm vol}\,B_{r}(0)} 
   \sum_{x \in \Lambda \cap B_{r}(a)} p(x) \; = \;
   \frac{1}{|\tilde{L}|} \, \int_{H} f(y) 
   \, {\rm d} \theta(y) $$
   uniformly in $a$. 
\end{theorem}
{\sc Proof}: The strategy will be to derive this more general result from
Proposition \ref{SchlottmannThm}. To this end, we approximate $f$ by a step
function $\psi$ on a (Riemann) admissible partition $\{U_1,\dots,U_n\}$ of
$W$, i.e.\ $W=\bigcup_{i=1}^n U_i$ with pairwise disjoint sets
$U_i\subset W$ that are all Riemann measurable.

{}Fix $\varepsilon>0$. By Lemma~\ref{step}, there is a step function 
$\psi$ on such an admissible partition (with suitable $n=n(\varepsilon)$)
of $W$, $\psi=\sum_{i=1}^n c_i \,\bs{1}^{}_{U_i}$, with 
$\|f-\psi\|_{\infty}<\varepsilon$. Choose a radius $R$ big enough so that 
we have, for all $r>R$, 
$$\left|\,\frac{1}{{\rm vol}\,B_r(0)}
   \left(\sum_{x\in\Lambda\cap B_r(a)} 1 \right) 
   - \, \frac{\theta(W)}{|\tilde{L}|}\,\right|
   \; < \; \varepsilon $$
and also
$$  \left|\,\frac{1}{{\rm vol}\,B_r(0)}
   \left(\sum_{x\in\Lambda(U_i)\cap B_r(a)} 1 \right) 
   - \,\frac{\theta(U_i)}{|\tilde{L}|}\,\right|
   \; < \; \frac{\varepsilon}{n}  $$
for all $1\le i\le n$, uniformly in $a$. Such a radius clearly exists. Then 
we have, since $p(x)=f(x^*)$, for all $r>R$ the following $3\varepsilon$-type 
argument,
\begin{eqnarray*}
\lefteqn{\left|\frac{1}{{\rm vol}\,B_r(0)}
         \left(\sum_{x\in\Lambda\cap B_r(a)} p(x)\right) 
         - \frac{1}{|\tilde{L}|}\,\int_W f(y) \, {\rm d}\theta(y)\right|} \\
   & \le & \frac{1}{{\rm vol}\,B_r(0)}
     \left(\sum_{x\in\Lambda\cap B_r(a)}
     \bigl| f(x^*) - \psi(x^*) \bigr| \right) \; + \;
     \frac{1}{|\tilde{L}|}\,\left|\int_W \bigl(f(y)-\psi(y)\bigr)
     \, {\rm d}\theta(y)\right| \; + \\
   & & \sum_{i=1}^n \, \left| \, \frac{c_i}{{\rm vol}\,B_r(0)}
     \left(\sum_{x\in\Lambda(U_i)\cap B_r(a)} 1 \right)
     - \frac{1}{|\tilde{L}|}\,\int_{U_i}\psi(y)\, {\rm d}\theta(y)\,\right|\\
   & < & \frac{\varepsilon}{{\rm vol}\,B_r(0)}
     \left(\sum_{x\in\Lambda\cap B_r(a)}\! 1\right) \; + \;
     \frac{\theta(W)}{|\tilde{L}|}\,\varepsilon \; + \\
   & & \sum_{i=1}^n \, |c_i| \, \left|\,\frac{1}{{\rm vol}\,B_r(0)}
     \left(\sum_{x\in\Lambda(U_i)\cap B_r(a)}\! 1\right)
     - \frac{\theta(U_i)}{|\tilde{L}|}\,\right| \\
   & < & \Bigl(2\,\frac{\theta(W)}{|\tilde{L}|} + 1\Bigr)
       \, \varepsilon \, + \,
       \sum_{i=1}^n \, \|\psi\|_{\infty} \, \frac{\varepsilon}{n} 
   \;\, < \;\, \Bigl(2\,\frac{\theta(W)}{|\tilde{L}|}
     + \|f\|_{\infty} + 2 \Bigr)\, \varepsilon\, ,
\end{eqnarray*}
from which the Theorem follows. \qed

\smallskip

{\sc Remark}: Weyl's Theorem also extends to all functions that are only
continuous on the compact set $W$, see the Remark following the proof of 
Lemma~\ref{step}.

\section{Self-similar densities on model sets} \label{model-dens}

Finally, we have collected all results that we need to construct
self-similar measures for model sets on their ``internal'' side,
and then, under certain circumstances, also invariant densities
on the model sets themselves.

\subsection{Self-similar systems}

An affine self-similar system of model sets consists of the the following 
data ({\bf SS1}--{\bf SS4}):
\begin{enumerate}
\item[\bf {SS1}] a cut and project scheme (\ref{diagram})
   whose internal space $H$, in addition to being an LCAG, is a complete 
   metric space with translation invariant metric $d$.
\item[{\bf{SS2}}] a family of regular model sets $\Lambda_i = \Lambda(W_i)$,
   $i=1,\ldots,n$, for this cut and project scheme, with each $W_i$ compact.
\item[{\bf{SS3}}] an invertible linear mapping $Q:\; \RR^m\rightarrow\RR^m$ which
   satisfies $Q(L)\subset L$, where $L$ is the projection $\pi(\tilde{L})$
   of the lattice $\tilde{L}$ in (\ref{diagram}). 
\item[{\bf{SS4}}] sets $F_{ij}$, $1\le i,j\le n$, some of which may be empty, of
   affine mappings
   $$C \; = \; C_a \;:\quad x \mapsto Qx+a \quad (a\in L)$$
   which map $\Lambda_j$ to $\Lambda_i$ and satisfy
\begin{equation}\label{eq7.1}
   \Lambda_i \; = \; \bigcup_{j=1}^{n} \;
    \bigcup_{C\in F_{ij}} C(\Lambda_j) \, , \quad 1\le i\le n\, .
\end{equation}
\end{enumerate}

The sets $F_{ij}$ may (and usually will) be infinite. Because all the
affine mappings involved have the same linear part, $Q$, each $F_{ij}$ is
parameterized by the translational parts $a\in L$. In the sequel, we will
thus mostly view $F_{ij}$ as a subset of $L$. In this case, we will, from
now on, use the notation $F_{ij}'$, i.e.\ $F_{ij}' := F_{ij}(0)$.

Such systems of model sets can arise quite naturally in the study of
self-similar tilings. \index{tiling} Each proto-tile is marked in some 
suitable way with a finite set of points, call them proto-points or 
{\em control points}. \index{control points}
This provides a marking of the tiling by points, and the sets $\Lambda_i$
are then taken to be the set of points that correspond to each class of
control points. In this case, the union in (\ref{eq7.1}) would typically
be disjoint, but in our study we definitely wish to include non-disjoint
unions as well.

The idea of this Section is to pass the self-similar system to the
internal side $H$ of the cut and project scheme, to apply our theory of
self-similar measures there, and finally to pull back the results to the
physical side, namely to the model sets $\Lambda_i$ themselves. We will see
that pulling back is not automatically possible, and we need to make various
types of assumptions to guarantee it. Still, these assumptions are not
unnatural, and they are actually met in many interesting cases. 

The situation for simple model sets $\Lambda = \Lambda(W)$ is just the  special
case of the general situation here, where $n = 1$. In this case, the matrix 
$\bs{s}$ that appears below is simply the unit matrix $(1)$.

We can directly lift $Q:\; L\rightarrow L$ to a group homomorphism
$\tilde{Q}:\;\tilde{L}\rightarrow\tilde{L}$, and then to a group homomorphism
$Q^*:\; L^*\rightarrow L^*$. We assume
\begin{itemize}
\item[\bf SS5] $Q^*$ is contractive with respect to the metric $d$.
\end{itemize}
In this case, since $L^*$ is dense in $H$, $Q^*$ extends to a continuous
contractive automorphism $$A:\quad H\longrightarrow H$$
with $A|^{}_{L^*}=Q^*$. Due to Lemma \ref{contraction-modulus}, the modulus 
$\alpha$ of $A$ with respect to the Haar measure $\theta$ on $H$ satisfies 
$\alpha>1$, see Section \ref{affine-maps} for details.

{}For each affine map $C_a:\; x\mapsto Qx+a$, $a\in L$, we define the 
affine mapping $C^*_a$ on $H$ by $y\mapsto Ay+a^*$. In this way, we arrive
at admissible families of contractions $F^*_{ij}$ on $H$ 
(see Section \ref{comp-fam}).
We let $\cF_{ij}$ be the closure of $F^*_{ij}$ in the space $C(H,H)$ of
continuous mappings on $H$. If we identify the mappings in $F_{ij}$ with 
their translational parts in $L$, then $F_{ij}^{\prime *}$ is viewed as a subset of
$L^*$ and $\cF_{ij}'$ is the closure of this in $H$ (see Section \ref{cfam}
where we did a similar thing). Let us summarize our notation in the
following diagram, where $*$ stands for the $*$-map and $'$ for the mapping
that links affine transformations with their translational parts.
$$\begin{array}{ccccccc}
  F^{}_{ij} & \stackrel{*}{\longleftrightarrow} & F^{*}_{ij} &
     \; \subset \; & \cF^{}_{ij}  & \; \subset \; & C(H,H) \\ &&&&\\
         ' \Big\updownarrow \hphantom{'} & & 
         ' \Big\updownarrow \hphantom{'} & &
         ' \Big\updownarrow \hphantom{'} & &
         \\ &&&&\\
  F^{\prime}_{ij} & \stackrel{*}{\longleftrightarrow} & F^{\prime *}_{ij} &
     \; \subset \; & \cF^{\prime}_{ij} & \; \subset \; & H
\end{array}$$

{}From (\ref{eq7.1}), we have, for all $1\le i\le n$,
\begin{equation}\label{eq7.2}
   \Lambda^*_i \;= \; \bigcup_{j=1}^{n}\;
    \bigcup_{C\in F_{ij}} C^*(\Lambda^*_j)
\end{equation} 
and taking closures gives us
\begin{equation}\label{eq7.3}
   W_i \; \supset \; \bigcup_{j=1}^{n}\;
    \bigcup_{C_{\vphantom{L}}^* \in F_{ij}^*} C^*(W_j) \, .
\end{equation}
Since the $W_i$ are compact and $C^*(W_j)=A W_j + a^* \subset W_i$ for
$C=C_a$, we see that the translational parts $F_{ij}'$ of the affine
mappings are bounded (with respect to $d$). Thus we have that
$\cF_{ij}'$ is compact in $H$  and $\cF_{ij}$ is compact in $C(H,H)$.

\begin{prop}\label{prop7.1} Under the above conditions, we have
\begin{equation}\label{eq7.4}
   W_i \; = \; \bigcup_{j=1}^{n}\;
      \bigcup_{D\in \cF_{ij}} D(W_j) \; , \quad i=1,\ldots,n\, .
\end{equation}
and $W_1 \times \ldots\times W_n$ is the attractor of $\cF$.
\end{prop}

To prove this result, we first establish
\begin{lemma}\label{lemma7.1}
Let $F$ be a relatively compact set of continuous mappings from $H$ to $H$.
Suppose that $U,V$ are compact subsets of $H$ such that $C(U)\subset V$
for all $C\in F$. Then, $D(U)\subset V$ for all $D\in\overline{F}$.
\end{lemma}

{\sc Proof}: Let $D\in\overline{F}$. Fix any $\varepsilon>0$ and let
$K=K(U,B_{\varepsilon}(0))$ be the set of all continuous mappings of $H$ to
itself that map $U$ inside $B_{\varepsilon}(0)$. This is an open neighbourhood
of $0$ in $C(H,H)$, so $D-C \in K$ for some $C\in F$. Thus
$$D(U)\;\subset\; C(U)+B_{\varepsilon}(0)\;\subset\;
  V+B_{\varepsilon}(0)\;\subset\; [V]^{}_{\varepsilon}\, . $$
This being true for all $\varepsilon>0$, we have $D(U)\subset V$.   \qed
\smallskip

{\sc Proof} of Proposition \ref{prop7.1}:
Consider (\ref{eq7.3}). Using Lemma \ref{lemma7.1}, we get
$$W_i\;\supset\;\bigcup_{j=1}^{n}\;\bigcup_{D\in\cF_{ij}}D(W_j)\,.$$
The right hand side is compact by (\ref{compactUnionMap}) and contains
$\Lambda^*_i$ by (\ref{eq7.2}), hence also $\overline{\Lambda^*_i}=W_i$. 
\qed \smallskip

{\sc Remark}: A solution to (\ref{eq7.4}) is guaranteed by the general
theory of contractions of Section \ref{mcms}. However, in the present
situation, we know more: the $W_i$ have non-empty interiors, since the
$\Lambda_i$ are model sets. In general, it is hard to know when such a
self-similar system of mappings leads to an attractor with non-empty
interior.

\subsection{Self-similar measures}

We now assume that each $\cF_{ij}$ is equipped with a positive regular
Borel measure $\sigma^{}_{ij}$, with $\sigma^{}_{ij}:=0$ if
$\cF_{ij}=\emptyset$ as before. As in Proposition \ref{aff-inv-meas},
we usually identify $\sigma^{}_{ij}$ with a positive Borel measure on $H$
that is supported on $\cF_{ij}'$. We set 
$s^{}_{ij}=\sigma^{}_{ij}(\cF_{ij}')=\sigma^{}_{ij}(H)$ and restate the
compatibility assumption {\bf CA} for the matrix 
$\bs{s}=(s^{}_{ij})$, namely that it has a positive 1-eigenvector:
\begin{itemize}
\item[\bf SS6] There is a positive vector $\bs{m}=(m^{}_1,\ldots,m^{}_n)^T$
  which satisfies $\bs{s}\bs{m}=\bs{m}$.
\end{itemize}

By Proposition \ref{aff-inv-meas}, we have the existence of a positive measure
$\omega^{}_1\otimes\ldots\otimes\omega^{}_n$, supported on
$W_1\times\ldots\times W_n$, which satisfies
\begin{equation}\label{eq7.5}
   \omega^{}_i \; = \; \sum_{j=1}^{n} \;\int_{\cF_{ij}}
      (D.\omega^{}_j )\, {\rm d}\sigma^{}_{ij}(D)\, ,
\end{equation}
with $\omega^{}_i(H)=m^{}_i$, for all $i=1,\ldots,n$, and which is explicitly 
given by the infinite product formula in Proposition \ref{aff-inv-meas}.
The next task is to convert (\ref{eq7.5}) into a statement about densities on
the model sets $\Lambda_1,\ldots,\Lambda_n$. Our basic assumption is
\begin{itemize}
\item[\bf SS7] Each $\omega^{}_i$ is {\em continuously representable}
  \index{continously representable} in the sense that, for $i=1,\ldots,n$,
  \begin{equation}\label{eq7.6}
       \omega^{}_i \; = \; g^{}_i\theta
  \end{equation}
  where $g_i\ge 0$ is a function on $H$ which is supported on $W_i$ and
  has the property that its restriction to $W_i$ is continuous on $W_i$.
\end{itemize}

Given (\ref{eq7.6}), we define the corresponding weights or
{\em densities} \index{densities}
$$p_i:\, L\rightarrow \RR \; , \quad p_i(x) \, = \, g_i(x^*)
  \; , \quad i=1,\ldots,n \, . $$
Since $g_i$ is supported on $W_i$, $p_i$ is supported on
$\Lambda_i=\{x\in L\mid x^*\in W_i\}$. Our assumption that each
$\Lambda_i$ is regular, i.e.\ that each $W_i$ is Riemann measurable,
allows us to apply Weyl's Theorem (Theorem \ref{weyl-thm}) to prove
the existence of the average density for each $p_i$:
\begin{equation}\label{eq7.7}
    \lim_{r\to\infty} \frac{1}{{\rm vol} (B_r(a))}
    \sum_{x\in L\cap B_r(a)} p_i(x) \; = \;
    \frac{1}{|\tilde{L}|} \int_H g_i \, {\rm d}\theta
    \; = \; \frac{m^{}_i}{|\tilde{L}|}
\end{equation}
where the convergence of the limit is uniform in $a$. 

{}For any affine mapping $D:\; x\mapsto Ax+a$ and any $h\in L^1(H)$,
we have $$D.(h\theta) \; = \; \alpha (D.h) \theta $$
as easily follows by applying both sides to a test function and then using
the formula ${\rm d}\theta(A^{-1}y) = \alpha\,{\rm d}\theta (y)$, see 
Proposition \ref{facts}(1). Plugging this into (\ref{eq7.5}) gives
$$g_i \; = \; \alpha \sum_{j=1}^{n}\int_{\cF_{ij}}
  (D.g_j)\,{\rm d}\sigma^{}_{ij}(D) $$
and then, for $i=1,\ldots,n$,
\begin{equation}\label{eq7.8}
  g_i(x) \; = \; \alpha \sum_{j=1}^{n}\int_{\cF_{ij}'}
      g_j\bigl(A^{-1}(x-u)\bigr)\,{\rm d}\sigma^{}_{ij}(u) \, .
\end{equation}

To pass this to the physical side, we need to be able to deal with the
integral. There are two situations in which we know how to do this, namely
when the $F_{ij}$ are finite and we basically use counting measures on the
$\cF_{ij}'$, and when the sets $\cF_{ij}'$ are Riemann measurable subsets
of $H$ and the measures $\sigma^{}_{ij}$ are basically restrictions of the 
Haar measure on $H$. Let us now discuss these cases.

\smallskip
\subsubsection{$F_{ij}$ finite, $\sigma^{}_{ij}$ counting measure}
If each set $F_{ij}$ is finite, then so is $F_{ij}^*$ and
$\cF_{ij}=\overline{F_{ij}^*}=F_{ij}^*$. The model sets $\Lambda_i$ are
linked by a finite collection of finite unions as follows,
\begin{equation}\label{eq7.9}
   \Lambda_i \; = \; \bigcup_{j=1}^{n} \bigcup_{k=1}^{N_{ij}}
   \,\bigl(Q \Lambda_j + a^{}_{ijk} \bigr)\, .
\end{equation}
We suppose that $\sigma^{}_{ij}$ is counting measure normalized to total
measure $s^{}_{ij}$ satisfying {\bf SS2}. Then
\begin{eqnarray}
  g_i & = & \alpha \,\sum_{j=1}^{n} \frac{s^{}_{ij}}{{\rm card} (F_{ij})}
            \sum_{a\in F_{ij}'} C_a^*.g_j   \label{eq7.10} \\
  p_i(x) & = & \alpha\, \sum_{j=1}^{n} \frac{s^{}_{ij}}{{\rm card} (F_{ij})}
            \sum_{a\in F_{ij}'} p_j\bigl(Q^{-1}(x-a)\bigr)\, .\label{eq7.11}
\end{eqnarray}

We do not know many conditions that guarantee the existence of the
representing functions $g_i$.

However, suppose that the unions in (\ref{eq7.2}) are {\em non-overlapping};
\index{IFS!non-overlapping} more precisely, we assume that
\begin{enumerate}
\item[{\bf NO1}] $F_{ij}$ is finite
\item[{\bf NO2}] $m^{}_i = \theta(W_i)$
\item[{\bf NO3}] $\sigma^{}_{ij}$ is counting measure scaled by $\alpha^{-1}$, i.e.\
      $s^{}_{ij}=\sigma^{}_{ij}(F_{ij})=\alpha^{-1}{\rm card}(F_{ij})$.
\item[{\bf NO4}] for each $i$, the sets $D(W_j)$ entering into the union in
      (\ref{eq7.4}) intersect at most on sets of measure $0$, i.e.\
      they are {\em just touching}. \index{IFS!just touching}
\end{enumerate}

Taking measures in (\ref{eq7.4}), we see how the compatibility condition
{\bf SS3} fits in:
$$m^{}_i \; = \; \sum_{j=1}^{n}\; \sum_{D\in \cF_{ij}} \alpha^{-1} m^{}_j
  \; = \; \alpha^{-1}\sum_{j=1}^{n} {\rm card} (F_{ij})\, m^{}_j
  \; = \; \sum_{j=1}^{n} s^{}_{ij} m^{}_j \, . $$
The self-similar measures $\omega^{}_i$ are easy to find, 
$\omega^{}_i = {\bf 1}^{}_{W_i} \theta$. In fact, the equation
$${\bf 1}^{}_{W_i} \; = \; \sum_{j=1}^{n}\; \sum_{D\in \cF_{ij}}
   {\bf 1}^{}_{D(W_j)} \qquad (a.e.) $$
(which is what non-overlapping means in (\ref{eq7.4})) is equivalent to
$${\bf 1}^{}_{W_i} \theta \; = \; \alpha^{-1}
  \sum_{j=1}^{n}\; \sum_{D\in \cF_{ij}} D.({\bf 1}^{}_{W_j} \theta)
  \; = \; \sum_{j=1}^{n}\;\int_{\cF_{ij}} D.({\bf 1}^{}_{W_j} \theta)
  \, {\rm d}\sigma^{}_{ij}(D)\, . $$
Consequently, the self-similar densities $p_i$ are simply
$$ p_i(x) \; = \;
   \begin{cases} 1 & \text{if $x\in\Lambda_i$} \\
                 0 & \text{otherwise}
   \end{cases} $$
and (\ref{eq7.11}), for $i=1,\ldots,n$, reads as
\begin{equation}\label{eq7.12}
  p_i(x) \; = \; \sum_{j=1}^{n}\;\sum_{a\in F_{ij}'}
      p_j\bigl(Q^{-1}(x-a)\bigr) \qquad (a.e.)
\end{equation}
where ($a.e.$) means that equality holds for $x\in L$, possibly up
to a set of density $0$.

Equation (\ref{eq7.12}) is the point set analogue of the situation in an
inflation tiling\index{inflation}\index{tiling!inflation} where the 
measure of the inflated tile is the sum of the
measures of the tiles into which it decomposes. Here, density replaces
measure and the ``tiling'' condition is effectively put on the attractor
in our assumptions in equations (\ref{eq7.4}).

Let us briefly mention that Lagarias and Wang \cite{LW} have recently 
begun an investigation
of multi-component point sets with self-similarities. Their paper
is taken from the discrete and combinatorial point of view,
but does address the issue of counting multiplicities due to
overlapping in the substitution process, and hence implicitly
the question of self-similar measures.

\smallskip
\subsubsection{$\cF_{ij}'$ Riemann measurable, $\sigma^{}_{ij}$ Haar measure}
By definition,
$$F_{ij}' \;\subset\;\{ a \in L\mid Q\Lambda_j + a \subset \Lambda_i \}\, ,$$
and so $F_{ij}^{\prime *}$ is a subset of
\begin{equation}\label{eq7.13}
   \cG_{ij}' \; := \; \{ b\in H\mid A W_j + b \subset W_i \} \, .
\end{equation}
Since $\cG_{ij}' = \bigcap_{w\in W_j} (W_i - A w)$ is closed (and compact), we
have $\cF_{ij}'\subset \cG_{ij}'$ for all $i,j$. Thus the $\cG_{ij}'$ give us 
an upper bound on the $\cF_{ij}'$ and, in any case,
$$F_{ij}' \; \subset \; \{ a\in L\mid a^*\in \cG_{ij}' \} \, .$$
The right hand side has the interesting property that, provided that 
$\cG_{ij}'$ is equal to the closure of its interior, it constitutes a model
set of the cut and project scheme (\ref{diagram}). If, furthermore, the
$\cG_{ij}'$ are Riemann measurable, then we have access to Weyl's Theorem
again. This suggests that we may use {\em all possible} self-similarities of a
given collection of model sets, replacing the $F_{ij}'$ by the sets
\begin{equation}\label{eq7.14}
   G_{ij}' \; = \; \{ a\in L\mid a^*\in\cG_{ij}'\}\, .
\end{equation}

Thus given regular model sets $\Lambda_i=\Lambda(W_i)$, $i=1,\ldots,n$,
and an inflation $Q$ satisfying {\bf SS1}--{\bf SS4} above,
we can define $\cG_{ij}'$ by (\ref{eq7.13}) and replace the given system 
$\cF$ of affine mappings by the new set $\cG$ with
\begin{equation}
   G_{ij} \; := \; \{ C: x \mapsto Qx + a \mid a \in L, a^* \in \cG_{ij}' \}\, .
\end{equation}

With this motivation in mind, we go back to our original setup with
{\bf SS1}--{\bf SS5}. We now assume in addition that
\begin{itemize}
\item[\bf SS8] each $\cF_{ij}'$ is Riemann measurable and, if
  $\theta(\cF_{ij}') > 0$, we have
  $$\sigma_{ij}^{} \; = \; \frac{s^{}_{ij}}{\theta(\cF_{ij}')}
    \,{\bf 1}^{}_{\cF_{ij}'}\, \theta \, ,$$
  i.e.\ $\sigma_{ij}^{}$ is Haar measure restricted to $\cF_{ij}'$ and
  normalized to total measure $s^{}_{ij}$ where $\bs{s}=(s^{}_{ij})$
  is an arbitrary non-negative matrix 
  satisfying {\bf SS2}. If $\theta(\cF_{ij}')=0$, 
  $\sigma^{}_{ij}$ is defined as the $0$-measure.
\end{itemize}

We now apply the results of Section \ref{mul-comp}. The $\sigma^{}_{ij}$
are absolutely continuous with respect to Haar measure $\theta$ since
$\sigma^{}_{ij}=h^{}_{ij} \theta$ where
$h^{}_{ij}=\frac{s^{}_{ij}}{\theta(\cF_{ij}')}{\bf 1}^{}_{\cF_{ij}'}$
and the density is $h_{ij}\in L^1(H)\cap L^{\infty}(H)$.

Assuming that the convolution (\ref{star}) converges
and using Proposition \ref{cont-inv-meas}, we obtain a family of
{\em continuous} non-negative functions $g^{}_1,\ldots,g^{}_n$
which satisfy the self-similarity equations
\begin{equation}\label{eq7.15}
  g^{}_i(x) \; = \; \alpha\,\sum_{j=1}^{n} w^{}_{ij} \,
   \int_{\cF_{ij}} g^{}_j \bigl( A^{-1} (x-u)\bigr)
   \, {\rm d}\theta(u)
\end{equation}
where $w^{}_{ij} := s^{}_{ij}/\theta(\cF_{ij}')$.
Applying Theorem \ref{weyl-thm} we obtain
\begin{theorem}\label{thm7.1}
   Let $\Lambda_1,\ldots,\Lambda_n$ be a self-similar system of model sets
   which satisfy the assumptions {\bf SS1}--{\bf SS6} and {\bf SS8}.
   Then there exist non-negative functions $p_i: L\rightarrow \RR$,
   supported on $\Lambda_i$, $i=1,\ldots,n$, with the following
   properties
\begin{eqnarray}\label{self-sim-den-eqn}
 m^{}_i & = & \lim_{r\to\infty}\, \frac{|\tilde{L}|}{{\rm vol}(B_r(a))}
       \sum_{x\in L\cap B_r(a)} p^{}_i(x) \, , \\
 p^{}_i(x) & = & 
       \lim_{r\to\infty} \,\frac{\alpha\,|\tilde{L}|}{{\rm vol}(B_r(a))}
       \;\sum_{j=1}^{n} w^{}_{ij} \sum_{x\in L\cap B_r(a)}
       p^{}_j\bigl(Q^{-1}(x-u)\bigr)\, ,
\end{eqnarray}
for all $x\in L$ and for all $i=1,\ldots,n$, where the limits are
uniform in $a\in\RR^m$. \qed
\end{theorem}

This may be compared with the similar formula that we derived in
\cite{BMmulti}. There, $\tilde{Q}$ was assumed to be an automorphism of
$\tilde{L}$, and the set of $Q$-affine mappings was the set of all
possible mappings. There, however, the scaling constants $\nu^{ij}$ 
(which are our $s^{}_{ij}$ here) had no general interpretation.

\smallskip
{\sc Remark}: For simplicity, the assumptions made in {\bf SS8} are actually 
a little stronger than necessary -- there is no general need to exclude the 
case that some
$\cF_{ij}'$ are singleton sets, but still carry a positive (point) measure.
Though this is then not absolutely continuous, it is `harmless' in the
convolution process because the convolution of a function with a unit point
mass only results in a shift of the function. We will meet this case in 
Section \ref{8.1.4} below.

\section{Concrete examples}\label{sec-examples}

A number of explicit examples have been presented earlier, in \cite{BM}
and \cite{BMmulti}. The case of the planar Penrose pattern
in its rhombic version is an example of a multi-component model set,
because the vertex points fall into 4 classes. This was described in
detail in \cite{BMmulti} and will not be repeated here. Instead, let us
look into a few other examples. First, we describe the silver mean chain 
\index{chain!silver mean} in one dimension and look at it in different ways, 
both as a single and as a multi-component model set. Next, we briefly describe
the Ammann-Beenker pattern \cite{AGS} \index{Ammann-Beenker pattern} in the 
plane, an eightfold symmetric relative of the Penrose pattern. This example 
appears also in other contributions to this volume.
Finally, we look at a more unusual example that involves the $3$-adic 
\index{3-adic numbers} numbers.

\subsection{The silver mean chain}

The silver mean chain is a $2$-sided sequence on the alphabet $\{a,b\}$. It can
generated by iterating the substitution
$$ a \mapsto aba \quad , \quad b \mapsto a $$
which, when starting from $a|a$, leads to the palindromic fixed point
\begin{equation}
 \dots   abaaabaabaabaaaba | abaaabaabaabaaaba \dots
\end{equation}
where $|$ simply marks the centre. With $a$ and $b$ interpreted as intervals 
($=$tiles) of length $\alpha:=1 + \sqrt 2$ and $1$ respectively, this gives 
rise to a tiling of the line, see \cite{HRB} for details of its structure.
The number $\alpha$ is called the {\em silver mean}. \index{silver mean}

Let $\Lambda_1$ (resp.\ $\Lambda_2$) denote the coordinates of the left  
end points of the tiles of type $a$ (resp.\ $b$), 
assuming that the initial block $a|a$ was centred
at $0$, i.e.\ its left end points are located at $-\alpha$ and $0$.  Then  
$\Lambda_1, \Lambda_2$, and $\Lambda := \Lambda_1 \cup \Lambda_2$ are 
model sets for the following cut and project scheme:
\begin{equation}\label{diagram2}
\begin{array}{ccccc}
   \RR & \stackrel{\pi}{\longleftarrow} &
   \RR\times \RR & \stackrel{\pi^{}_\RR}{\longrightarrow} & \RR \\
   \cup & & \cup & & \cup \\
   \ZZ[\sqrt{2}\,] & \longleftarrow & \vphantom{\bigcap^N}
   \tilde{L} & \longrightarrow &  \ZZ[\sqrt{2}\,] \\    
\end{array}
\end{equation}
where $\ZZ[\sqrt{2}\,]=\ZZ\oplus\ZZ\sqrt{2}$ is the ring of integers in the
quadratic field $\QQ(\sqrt{2}\,)$, the $^*$-map from $\ZZ[\sqrt{2}\,]$ to 
$\ZZ[\sqrt{2}\,]$ is the algebraic conjugation defined by 
$\sqrt 2 \mapsto -\sqrt 2$, and the lattice is
$\tilde{L} := \{(x, x^*) \, | \, x \in \ZZ[\sqrt{2}\,] \}$.
Explicitly, we obtain
\begin{eqnarray}\label{silver-mean}
\Lambda_1  &=& \{ x \in \ZZ[\sqrt{2}\,] \mid x^* \in W_{1} \} \nonumber \\
\Lambda_2  &=& \{ x \in \ZZ[\sqrt{2}\,] \mid x^* \in W_{2} \}           \\
\Lambda_{\hphantom{1}} &=& 
         \{ x \in \ZZ[\sqrt{2}\,] \mid x^* \in W_{\hphantom{1}} \} \nonumber
\end{eqnarray}
where the corresponding windows are intervals, namely
\begin{equation} \label{windows}
     W_1 \; := \; \Bigl[\frac{1}{\sqrt 2} -1,\frac{1}{\sqrt 2} \,\Bigr]\, , \; 
     W_2 \; := \; \Bigl[-\frac{1}{\sqrt 2},\frac{1}{\sqrt 2}-1\,\Bigr]\, , \;
     W   \; := \; \Bigl[-\frac{1}{\sqrt 2},\frac{1}{\sqrt 2}\,\Bigr]\, ,
\end{equation}
where $W=W_1\cup W_2$ and $W_1\cap W_2 = \{\alpha^*/\sqrt{2}\,\}$.

This can be verified in the following way. Let $Q$ be the linear mapping which 
is scalar multiplication by $\alpha = 1 +\sqrt2$ and let $A$ be its conjugate 
map, which is scalar multiplication by $\alpha^* = -\alpha^{-1} = 1-\sqrt{2}$. 
With the explicit coordinatization given above, the substitution rules say that
\begin{equation}\label{silver-mean-Substrules}
\begin{split}
  \Lambda_1  \;&=\;  Q\Lambda_1  \cup Q\Lambda_2  \cup (Q\Lambda_1 + \alpha +1)\\
  \Lambda_2  \;&=\;  Q\Lambda_1 + \alpha
\end{split}
\end{equation}
so that we now have a system of $Q$-inflations which, in the notation of 
{\bf SS4}, is defined by
\begin{equation}\label{substitution-matrix}
  F' \; = \; (F'_{ij}) \; = \; 
  \begin{pmatrix}
    \{0, \alpha +1\} & \{0 \} \\
    \{\alpha \} & \emptyset
  \end{pmatrix}. 
\end{equation}
The corresponding contractions satisfy:
\begin{equation}\label{silver-mean-SubstrulesInternal}
  \begin{split}
    W_1 \;&=\; AW_1 \cup AW_2  \cup (AW_1 + \alpha^* +1)\\
    W_2 \;&=\; AW_1 + \alpha^*
  \end{split}
\end{equation}
which shows that $W_1 \times W_2$ is the attractor for the system of $A$-affine
contractions given by $F^{\prime *}$.

Since the generators of $\Lambda$ (which correspond to the $a$-tiles with 
coordinates $0$ and $-\alpha$) are mapped into  $W_1$ by $^*$, all subsequent 
points generated from them are $^*$-mapped into $W$. Thus the model set 
$\Lambda(W)$ assuredly contains our set $\Lambda$. On the other hand, 
it is easy to see that the minimum separation between the points of $\Lambda(W)$
is $1$, and no point can be added to $\Lambda$ without violating this. In short, 
$\Lambda = \Lambda(W), \Lambda_1 = \Lambda(W_1), \Lambda_2 = \Lambda(W_2) $.

We now examine this situation in four different ways, namely as
single component and multi-component case, and each then with minimal
and maximal families of affine contractions. By this we mean either
the case when the window system is minimally generated by affine contractions or
when we use the entire set of affine contractions available.

We content ourselves with a few remarks about each case and one figure illustrating
the continuous case. We indicate the appropriate Sections of the paper as we go 
along and freely use the notation from these Sections. Note that 
$\alpha = 1 + \sqrt 2$ as we have defined it in this
Section is the modulus of the contraction $A$, in keeping with previous notation.

\smallskip
\subsubsection{The single model set $\Lambda = \Lambda(W)$ with $F$ minimal} 
(see Section 3.2) 

\noindent The contractivity factor is $\alpha^*$. Since $|\alpha^*| < 1/2$, 
we need at least three affine mappings to get the full window as the attractor --
we would end up with a Cantor subset of $W$ if we would start with only two.
One possible choice is
$$ W \; = \; (AW + \alpha^*) \cup AW \cup (AW - \alpha^*) $$
and our family of mappings is then 
$F^{}_H = F^{}_{\RR} = \{ \alpha^*, 0, -\alpha^* \}$.
We take the simplest of all probability measures on $F^{}_{\RR}$, i.e.\
counting measure, 
$$
   \nu \; = \; \mbox{$\frac{1}{3}$}
            (\delta_{\alpha^*} + \delta^{}_0 + \delta_{-\alpha^*}) \, ,
$$
so that this is an example of a finite IFS.
The corresponding invariant measure on $W$ and its Fourier transform  are given by 
\begin{equation}
  \begin{split}
     \omega   \;&=\; \overset{\infty}{\underset{\ell=1}{\conv}}
                  \, \mbox{$\frac{1}{3}$}
                  \bigl(\delta_{(\alpha^*)^{\ell}}+\delta^{}_0
                        +\delta_{-(\alpha^*)^{\ell}}\bigr)  \\
     \hat{\omega} \;&=\; \overset{\infty}{\underset{\ell=1}{\mbox{$\prod$}}}
            \, \mbox{$\frac{1}{3}$}
            \bigl( 1 + 2 \cos (2 \pi i (\alpha^*)^{\ell} k) \bigr) \, .
  \end{split}
\end{equation}

This measure is similar to those studied in the context of the binary
addressing problem, see \cite{Sol} and references therein. Although most
of them are absolutely continuous if the IFS covers the full interval
(which is does here), exceptions emerge, see \cite[Thm.\ 4]{BG}, 
if the scaling factor of the IFS 
is the inverse of a Pisot-Vijayaraghavan number (which is the case here, too).
These exceptional self-similar measures will then be purely singular
continuous. Note that this might be very difficult to detect 
if one is not aware of it -- the fractal dimension of such a measure
can be tantalizingly close to 1, see \cite[Sec.\ 8]{Lal}. 
In any case, we cannot pull back such a 
measure to the physical side, and thus do not gain much insight into
the structure of our model set from it.

% This measure was studied in \cite{JW}. Since $|\alpha^*| = \sqrt 2 -1 < \frac{1}{2}$,
% the measure is singular with respect to Lebesgue measure on $\RR$. For other
% results in this direction, we refer to \cite{Sol}. Unfortunately, our theory does 
% not allow us to gain any special information on $\Lambda$ in this case.

\smallskip
\subsubsection{The single model set $\Lambda = \Lambda(W)$ with $F$ maximal} 
(see Section 3.3) 

\begin{figure}
\psfig{file=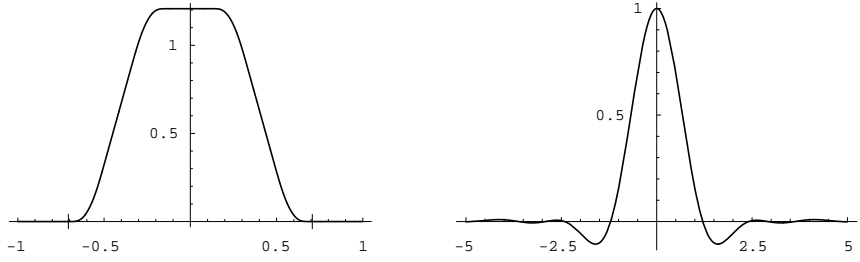,height=0.18\textheight}
\caption{Invariant density (left) and its Fourier transform (right) for the 
silver mean chain. The support of the $C^{\infty}$-density $g$ of
(\ref{ex-1.2}) is the window $W$ of (\ref{windows}), indicated by extra markers.
\label{SMfig1}}
\end{figure}

\noindent Here, we start with the observation
$$ W \; = \; \bigcup_{u \in [\alpha^*, -\alpha^*]} AW + u  \, . $$
We now use the complete set of available self-similarities 
$\cF' = [\alpha^*, -\alpha^*]$. For our pre-assigned probability 
measure $\nu$, we choose
$$
  \nu \; = \;  \frac{\bs{1}^{}_{[\alpha^*, -\alpha^*]}}{2|\alpha^*|} \, \theta
      \; = \;  \frac{\bs{1}^{}_{\cF'}}{2|\alpha^*|} \, \theta
$$
where $\theta$ is Lebesgue measure on $\RR$. Then we are the situation of 
Proposition \ref{single-conv}:
\begin{equation} \label{ex-1.2}
  \begin{split}
    g  \;&=\;  \overset{\infty}{\underset{\ell=0}{\conv}}  \left( 
    \frac{\bs{1}^{}_{|\alpha^*|^{\ell}\cF'}}
         {2\,|\alpha^*|^{\ell+1}} \right)   \\
    \hat{g}(k)  \;&=\;  
    \overset{\infty}{\underset{\ell=1}{\mbox{$\prod$}}} \left( 
           \frac{\sin(2 \pi (\alpha^*)^{\ell} k)}{2 \pi(\alpha^*)^{\ell} k}   
    \right) \, . 
  \end{split}
\end{equation}
This is illustrated in Figure \ref{SMfig1}.

\smallskip
\subsubsection{The system of model sets $\{\Lambda_1,\Lambda_2\}$ with
windows $\{W_1,W_2\}$, with $F$ minimal} 
(see Section 7.2, non-overlapping case)  \label{8.1.3}

We consider $\Lambda = \Lambda_1 \cup \Lambda_2$ as a
multi-component model set through (\ref{silver-mean-Substrules}) and 
(\ref{substitution-matrix}) with the windows
$W_1$ and $W_2$ given by (\ref{windows}) and related by 
(\ref{silver-mean-SubstrulesInternal}). The non-overlapping
conditions {\bf NO1} and {\bf NO4} hold. To obtain {\bf NO2} and {\bf NO3},
we define $\bs{m} = (1, |\alpha^*|)^T$ and (since $|\alpha^*| = \alpha - 2
= 1/\alpha$)
\begin{equation}
   \bs{s} \; = \; \begin{pmatrix}
                     2/\alpha &  1/\alpha\\
                     1/\alpha &  0
                  \end{pmatrix} ,
\end{equation}
which is clearly primitive.
The self-similar measures are $\omega_1 = \bs{1}^{}_{W_1} \theta$ and 
$\omega_2 = \bs{1}^{}_{W_2}\theta$, and the corresponding densities on the 
physical side are $p_1$ and $p_2$, which are the characteristic functions 
(defined on $\ZZ[\sqrt{2}\,]$) of $\Lambda_1$ and $\Lambda_2$, respectively.

\smallskip
\subsubsection{The system of model sets $\{\Lambda_1,\Lambda_2\}$ with
windows $\{W_1,W_2\}$, with $F$ maximal} 
(see Section 7.2, Riemann measurable case) \label{8.1.4}

We continue with the multi-component picture of Section \ref{8.1.3}, but now we 
use all available $\alpha$-affine self-similarities. These are easily computed 
from (\ref{windows}) and (\ref{silver-mean-SubstrulesInternal}), resulting in
$$ \cF' \; = \; \begin{pmatrix}
         [0, 1+\alpha^*] & [\alpha^*, -\alpha^*] \\
            \{\alpha^*\} & \emptyset
                \end{pmatrix} .
$$
The sets appearing here are all measurable. Let us fix the same vector 
$\bs m$ and matrix $\bs s$ as in Section \ref{8.1.3}. If we use Lebesgue measure 
$\theta$ on $\RR$, we obtain {\bf SS7} in the form
$$ \bs{\sigma} \; = \; \begin{pmatrix}
         (1 - \alpha^*)\bs{1}^{}_{[0,1+\alpha^*]} \, \theta  &   
         \frac{1}{2} \bs{1}^{}_{[\alpha^*,-\alpha^*]} \, \theta \\
         -\alpha^*\delta^{}_{\alpha^*} & 0
                  \end{pmatrix}.
$$
Using Proposition \ref{aff-inv-meas}, we have the self-similar measure
\begin{equation} \label{bigConv}
\omega \, = \, \overset{\infty}{\underset{\ell=0}{\conv}} \; \alpha^{\ell}
    \begin{pmatrix}
       (1 - \alpha^*)\bs{1}^{}_{(\alpha^*)^{\ell}[0,1+\alpha^*]} \, \theta  & 
      \frac{1}{2} \bs{1}^{}_{(\alpha^*)^{\ell}[\alpha^*,-\alpha^*]} \, \theta \, \\
      -\alpha^*\delta^{}_{(\alpha^*)^{\ell+1}} & 0
    \end{pmatrix}     \conv
    \begin{pmatrix} \delta^{}_0 \\ |\alpha^*|\, \delta^{}_0  \end{pmatrix} .
\end{equation}

The solution here is mildly different from the one appearing  in 
${\bf SS8}$ because of the appearance of the point measure, compare the
Remark following Theorem \ref{thm7.1}.
However, the convolution of a delta and a function 
simply translates the function. Furthermore, the convolution of two functions 
$f$ and $g$ is point-wise bounded by $||f||_\infty ||g||_1$.
Thus the partial convolution products in  $\omega^{(n)}$ of $\omega$ in 
(\ref{bigConv}), for all $n \ge 1$, are {\em functions} 
$$\begin{pmatrix} f^{(n)}_{11}  & f^{(n)}_{12} \\
                  f^{(n)}_{21}  & f^{(n)}_{22} \end{pmatrix}
$$
where the $f^{(n)}_{ij} $ are supported on $W_j$, are uniformly bounded,  
and are increasingly differentiable as $n \to\infty$. Alternatively, one
could also work with the square of the averaging operator here, which would
match ${\bf SS8}$ from the beginning.

The resulting self-similar measure is represented by two $C^{\infty}$ 
functions, $g^{}_1$ and $g^{}_2$. They are of the kind shown on the left 
of Figure \ref{SMfig1}, but now supported on $W_1$ and $W_2$, respectively,
and with total mass $1$ and $|\alpha^*|$, in accordance with $\bs{m}$.

\subsection{The Ammann-Beenker model set}

This relative of the rhombic Penrose tiling is usually described as a model
set obtained from the primitive cubic lattice $\ZZ^4$ in 4-space, see
\cite{BJ} and references therein for details. Here, we prefer the number
theoretic approach given in \cite{HRB}, which is more compatible with
the above description of the silver mean chain. With $\xi:= e^{2 \pi i/8}$,
we use the following cut and project scheme
\begin{equation}\label{diagram3}
\begin{array}{ccccc}
   \RR^2 & \stackrel{\pi}{\longleftarrow} &
   \RR^2\times \RR^2 & \stackrel{\pi^{}_{\RR^2}}{\longrightarrow} & \RR^2 \\
   \cup & & \cup & & \cup \\
   \ZZ[\xi] & \longleftarrow & \vphantom{\bigcap^N}
   \tilde{L} & \longrightarrow &  \ZZ[\xi] \\    
\end{array}
\end{equation}
where $\ZZ[\xi]$ is the ring of cyclotomic integers generated by the
primitive 8th roots of unity. It is the maximal order in the cyclotomic
field $\QQ(\xi)$. Then, the lattice is
$$\tilde{L} \; = \; \{(x,x^*)\mid x\in\ZZ[\xi] \}$$
where the $*$-map is given by algebraic conjugation $\xi\mapsto\xi^3$.

\begin{figure}
\mbox{\psfig{file=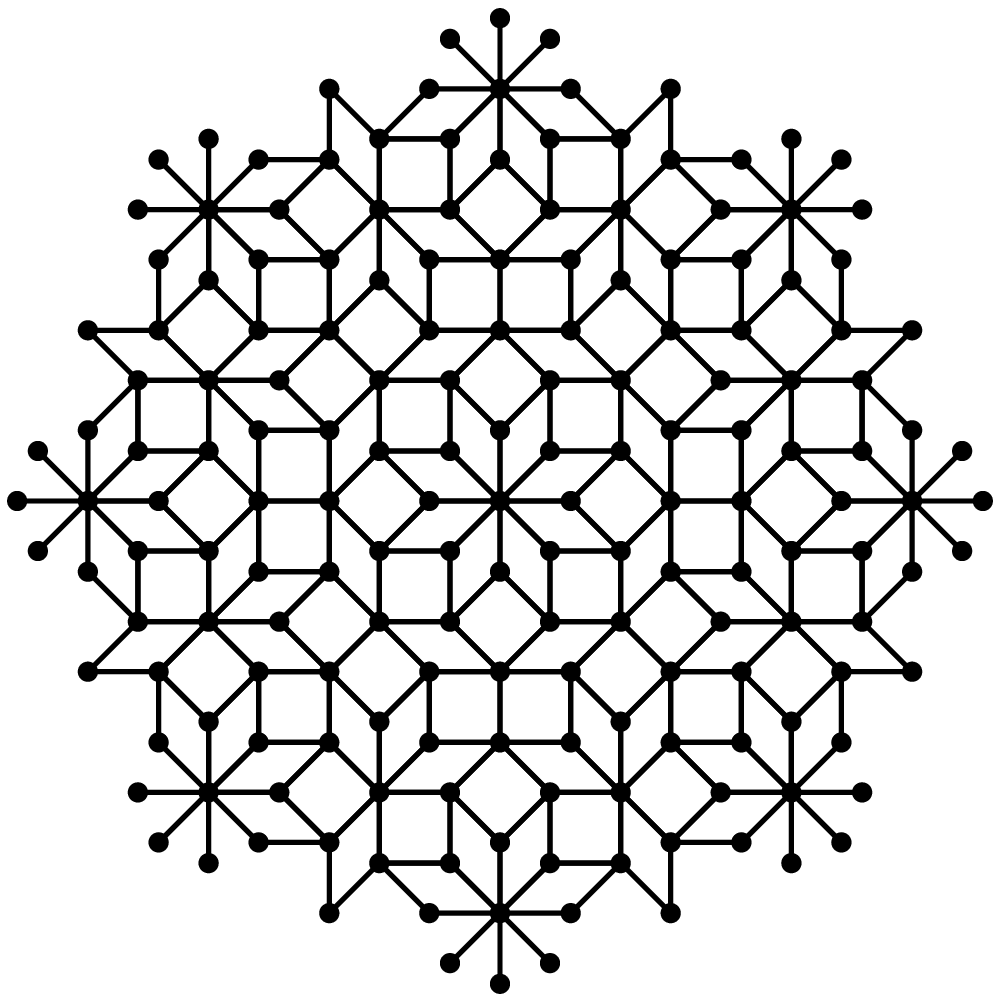, height=0.33\textheight}\hspace*{0.05\textwidth}
\psfig{file=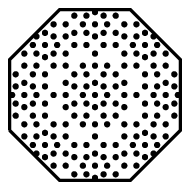,height=0.35\textheight}}
\caption{The Ammann-Beenker model set, seen as the vertex set of a tiling
with squares and rhombi. The window in internal space is a regular octagon.
It is shown to the right with the $*$-images of the vertex points from the
patch to the left.      \label{ABfig}}
\end{figure}

In this setting, the window $W$ of the Ammann-Beenker model set
$\Lambda = \Lambda(W)$ is simply a regular octagon of edge length 1,
centred at the origin. Its area is $2\alpha$, with $\alpha=1+\sqrt{2}$ 
as above. With this choice, the model set $\Lambda$ is both 
{\em regular}\index{model set!regular} ($W$ is a polytope, hence Riemann 
measurable) and {\em generic}\index{model set!generic}
($L=\ZZ[\xi]$ does not intersect $\partial W$). 
A symmetric patch and its lift to internal space 
is shown in Figure \ref{ABfig}.

Let us now consider $\Lambda=\Lambda(W)$ as a single model set and let us determine
the invariant density on $W$ that results from the set of {\em all}
self-similarities of the form $\Lambda \mapsto \alpha\Lambda + v\subset\Lambda$
(with $v\in\ZZ[\xi]$). So, $Q$ is again multiplication by $\alpha$, and $A$ 
multiplication by $\alpha^*$, a contraction. It then follows that
$$\cF' \; = \; \{ u \mid \alpha^* W + u \subset W \} 
       \; = \; (2-\sqrt{2}\,) W  \, ,$$
i.e.\ $\cF' \subset\RR^2$ is another octagon centred at the origin, but with 
reduced edge 
length $|\alpha^*|\sqrt{2} = 2-\sqrt{2}$. Consequently, $\cF'$ has area 
$4 |\alpha^*|$. For the a priori probability measure on $\cF'$, we choose 
$$\nu \; = \; \frac{\bs{1}^{}_{\cF'}}{4 |\alpha^*|}\,\theta$$
where $\theta$ is now Lebesgue measure on $\RR^2$. Proposition
\ref{single-conv} then gives
\begin{equation} \label{ab-dens-formula}
   g  \; = \;  \overset{\infty}{\underset{\ell=0}{\conv}} 
    \left( \frac{\bs{1}^{}_{|\alpha^*|^{\ell} \cF'}}
         {4\,|\alpha^*|^{2\ell+1}} \right) \, .
\end{equation}

\begin{figure}
\vspace*{-8mm}
\psfig{file=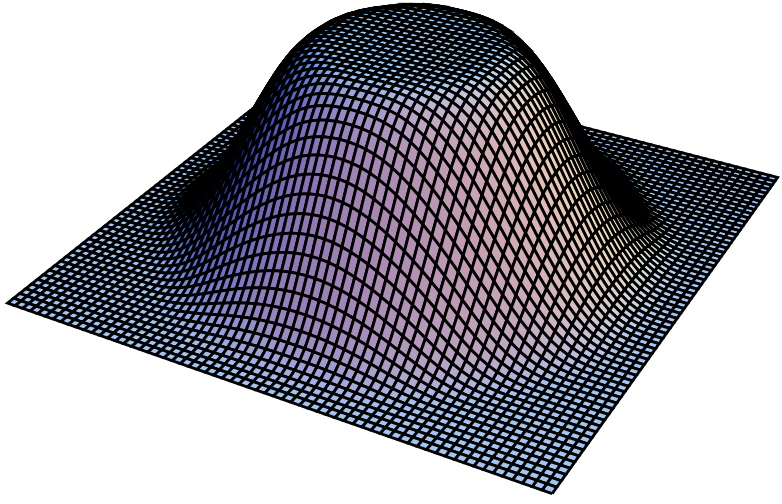,height=0.4\textheight}
\caption{The self-similar density $g$ for the Ammann-Beenker model set,
according to Eq.\ (\ref{ab-dens-formula}).
The support of the $C^{\infty}$-function $g$ is a regular octagon,
the window $W$ in the upper right corner of Figure \ref{ABfig}.
\label{ABdens}}
\end{figure}

The invariant density is a $C^{\infty}$-function supported on the window $W$,
see Figure \ref{ABdens}. The deviations from circular symmetry are rather faint.
However, in contrast to the Penrose case investigated in \cite{BMmulti}, it
has a central plateau and then rolls off smoothly towards the boundary of $W$.
The same phenomenon is actually visible for the invariant density of the
silver mean chain in Figure \ref{SMfig1}, in contrast to the one for the
Fibonacci chain in \cite{BM}. It is due to the larger absolute value of
$\alpha$ in comparison to the golden ratio $\tau = (1+\sqrt{5}\,)/2$ which
appears there. This is interesting in relation to the rather widespread 
experimental finding that ``real world'' quasicrystals are to be described
by window functions with a smooth roll-off. Although this is usually explained
as a random tiling effect, our above examples show that other mechanisms are
possible as well, and the attractive feature is then that they result from
some residual (or statistical) inflation\index{inflation symmetry} symmetry.

\subsection{A $3$-adic example}

In this final example, which will be fully developed along with some other 
$p$-adic examples in \cite{HM}, we indicate what invariant measures can look 
like in a very different situation. The result is quite surprising.

This time, we begin with the $2$-sided chain on the ternary alphabet
$\{a,b,c\}$, generated by the substitution rule
$$
a \mapsto ab \quad , \quad b \mapsto abc \quad , \quad c \mapsto abcc \, .
$$
Starting from $c|a$, iteration leads to the 2-sided fixed point
\begin{equation}
  \dots babcabccabccababcabccabcc | ababcababcabccababcababc \dots
\end{equation}
With $a,b,c$ assigned intervals of length $1,2,3$ respectively, which is
the natural geometric realization here, this again gives 
rise to a tiling of the line. This system was studied in \cite{BMS} where it was 
shown that, when it is coordinatized, the resulting sets of points 
$\Lambda_1, \Lambda_2, \Lambda_3$ are model sets based on the $3$-adic integers
$\hat{\ZZ}_3$. The coordinatization starting at $0$ and using {\em right} 
end points is
\begin{eqnarray}\label{3-adic-sets}
   \Lambda_1 &=& \mbox{$\bigcup_{k=2}^\infty $} \Bigl(
             \bigl( \mbox{$\sum_{i=0}^{k-2} $} \, 3^i \bigr) + 3^k\ZZ \Bigr)
             \nonumber \\
   \Lambda_2 &=& \mbox{$\bigcup_{k=2}^\infty $} \Bigl(
             \bigl(2 + \mbox{$ \sum_{i=0}^{k-2} $} \, 3^i \bigr) + 3^k\ZZ \Bigr)\\
   \Lambda_3 &=& 3^2\ZZ \cup \left(\mbox{$\bigcup_{k=3}^\infty $} \Bigl(
              \bigl(- \mbox{$\sum_{i=1}^{k-2} $}\, 3^i \bigr) 
              + 3^k\ZZ \Bigr) \right) \, .
              \nonumber
\end{eqnarray}

The corresponding model sets are 
$$ \Lambda_i = \{ x \in \ZZ \, | \, x \in W_i \}$$
where the windows are given by
\begin{eqnarray}\label{3-adic-windows}
 W_1 &=& \mbox{$ \bigcup_{k=2}^\infty $}\left( \bigl(
         \mbox{$\sum_{i=0}^{k-2} $}\, 3^i \bigr) + 3^k \hat{\ZZ}_3 \right)
         \nonumber \\
 W_2 &=& \mbox{$\bigcup_{k=2}^\infty $} \left( \bigl(
         2 + \mbox{$\sum_{i=0}^{k-2} $}\, 3^i \bigr)
         + 3^k \hat{\ZZ}_3 \right) \\
 W_3 &=& 3^2\hat{\ZZ}_3 \cup \left(\mbox{$\bigcup_{k=3}^\infty $} \Bigl( 
         \bigl(-\mbox{$\sum_{i=1}^{k-2} $} \, 3^i \bigr)
         + 3^k \hat{\ZZ}_3 \Bigr) \right) \, . \nonumber
\end{eqnarray}

We consider this as a multi-component model set and the basic inflationary 
maps as affine mappings of the form $x \mapsto 3x + u$, $ u \in \ZZ$. These 
maps transfer over to the internal side without any symbolic change, but 
there they are contractions with respect to the standard
$3$-adic metric topology.  The most interesting case seems to be when we 
allow all possible mappings of this type between the various windows 
(this corresponds to the case of the silver mean example
in Section \ref{8.1.4}). 
These are quite straightforward to work out and give rise to a 
$3\times 3$ matrix $\cF$. The invariant measures on  $\hat{\ZZ}_3$ are 
then given by the matrix equation
\begin{equation}
     \omega \; = \;  \left(\,\overset{\infty}{\underset{\ell=0}{\conv}} \,
          Y^{(\ell)}\right) * \, \delta^{\bs m}
\end{equation}
where 
\begin{equation}
    Y^{(\ell)} \; := \; \left( \frac{s^{}_{ij}}{\theta(3^{\ell} \cF_{ij})}\, 
            \bs{1}^{}_{3^{\ell} \cF_{ij}} \right)_{1\le i,j \le 3} \, ,
\end{equation}
$\theta$ is the Haar measure on $\hat{\ZZ}_3$ normalized to total
measure $1$, and $\bs{s} \bs m = \bs m$. 

Remarkably, this convolution can be explicitly computed and yields
\begin{equation}
   \omega \; = \; 3^2 \begin{pmatrix}
          \,\bs{1}^{}_{3^2\hat{\ZZ}_3 +1}  \, m^{}_1\, \theta \, \\
          \,\bs{1}^{}_{3^2\hat{\ZZ}_3 +3}  \, m^{}_2\, \theta \, \\
          \,\bs{1}^{}_{3^2\hat{\ZZ}_3}\quad\, m^{}_3\, \theta \,
                      \end{pmatrix} . 
\end{equation}
In particular, the self-similar measures are absolutely continuous and give 
rise to the following self-similar densities on our three model sets
$\Lambda_{\{1,2,3\}}$ in $\ZZ$:
$$
   p^{}_{\{1,2,3\}}(\ell) \, = \, \begin{cases} 
              9\, m^{}_{\{1,2,3\}} & \text{if $\ell\equiv \{1,3,0\}$ mod $9$}, \\
              0 & \text{otherwise}. \end{cases}
$$

Thus, the self-similar densities are {\em periodic} although the point sets 
$\Lambda_i$ are aperiodic. In fact, the support of the densities consists 
of three different cosets of  $9 \ZZ$, each of which lies entirely inside 
one of the three point sets $\Lambda_i$, $i=1,2,3$. We do not know many other 
mechanisms that result in periodic\index{periodic state} states, although it 
is not clear whether this inflation induced periodicity could be interpreted 
physically.

\renewcommand{\thesection}{\Alph{section}}
\setcounter{section}{1}

\section*{Appendix}

{}For the sake of completeness, we collect a number of results in this
Appendix. Let us first show that real (complex) Lipschitz functions are 
dense in $C(X,\RR)$ (in $C(X,\CC)$) if $X$ a compact metric space. We
first need the following result which is essentially stated in 
\cite[Prop.\ IX.2.2.3]{Bourbaki}.

\begin{lemma}\label{sepa}
Let $X$ be an arbitrary metric space.
{}For every non-empty $A\subset X$, the function $\phi\!:x\mapsto d(x,A)$ 
is in ${\rm Lip}(\le\!1,X,\RR)$. The Lipschitz constant is $r^{}_{\phi}=0$
iff $\phi\equiv 0$ on $X$, and $r^{}_{\phi}=1$ otherwise.
\end{lemma}
{\sc Proof}: Recall that $d(x,A)=\inf\{d(x,z)\mid z\in A\}$.
Let $x,y\in X$. {}Fix $\varepsilon >0$ and choose $z\in A$ such that 
$d(y,z)\le d(y,A) + \varepsilon$. Then
\begin{eqnarray*}
  d(x,A) - d(y,A) & \le & d(x,A) - d(y,z) + \varepsilon
      \; \le \; d(x,z) - d(y,z) + \varepsilon \\
  & \le & d(x,y) + \varepsilon \, .
\end{eqnarray*}
Now, since $\varepsilon >0$ was arbitrary and the argument is
essentially symmetric in $x$ and $y$, we can conclude
$$ |\phi(x)-\phi(y)| \; = \; |d(x,A) - d(y,A)| \; \le \; d(x,y) \, . $$
This means that $\phi:x\mapsto d(x,A)$ is Lipschitz with constant 
$r^{}_{\phi}\le 1$.

If $\phi\equiv 0$, we trivially have $r^{}_{\phi}=0$. Otherwise, there
exists an $x\in X\setminus A$ with $\phi(x)=d(x,A)>0$. We know
from above that $r^{}_{\phi}\le 1$. Fix $\varepsilon>0$ and $y\in A$ with 
$d(x,y)\le d(x,A)+\varepsilon$. Then, we have $d(y,A)=0$ and get
$$ 0 \; < \; d(x,y) \; \le \; d(x,A) + \varepsilon \; = \;
   d(x,A) - d(y,A) + \varepsilon \; \le \; 
   r^{}_{\phi} \, d(x,y) + \varepsilon \, . $$
Since $\varepsilon >0$ was arbitrary, this is only possible with
$r^{}_{\phi}\ge 1$, hence $r^{}_{\phi}= 1$.  \qed

\begin{lemma}\label{lip-real} If $X$ is a compact metric space,
the real $\,(\!$complex{\rm )} Lipschitz functions are dense in $C(X,\RR)$
$($in $C(X,\CC)$$)$.
\end{lemma}
{\sc Proof}:
This is a straight-forward application of the Stone-Weierstra{\ss} theorems.
The real Lipschitz functions form a subalgebra
of $C(X,\RR)$ under pointwise addition and multiplication because
$r^{}_{\phi+\psi} \le r^{}_{\phi} + r^{}_{\psi}$ and
$r^{}_{\phi\psi}  \le \|\phi\|^{}_{\infty} r^{}_{\psi} + 
\|\psi\|^{}_{\infty} r^{}_{\phi}$.
Furthermore, the constant functions are Lipschitz, and the separation 
property follows from Lemma~\ref{sepa}. So, the real variant, see 
\cite[Thm.\ IV.9]{RS}, gives the one claim, while the complex variant,
\cite[Thm.\ IV.10]{RS}, gives the other. In the latter case,  
the additional requirement that the algebra of complex Lipschitz functions 
is closed under complex conjugation is obvious.  \qed

\smallskip

We also need a number of convergence results,
which we collect here. The following standard result is a special case 
of \cite[Thm.\ 14.22 and Cor.\ 14.23]{Q} or
\cite[Thms.\ 7.14 and 7.15]{Kelley}, see also \cite[A 14.8]{Q}.
\begin{lemma}\label{point-compact}
Let $X$ be a locally compact space and $\{g_n\}$ an equi-continuous 
sequence of functions in $C(X,\CC)$ which converges pointwise.
Then the limit is a continuous function, $g\in C(X,\CC)$, and the
convergence is uniform on compact subsets of $X$ 
$\,(\/$compact convergence $\!)$.   \qed
\end{lemma}

\begin{lemma}\label{equi-cont}
Let $W$ be a compact subset of an LCAG $H$, with dual group $\widehat{H}$.
Let $(\mu_n)_{n\in\NN}$ be a vaguely convergent sequence of probability 
measures in $\cP(W)$, with limit $\mu$. If $K\subset\widehat{H}$ is compact, 
the family $\{\widehat{\mu}_n\}$ of functions from $C(H,\CC)$ is 
equi-continuous $\,(\/$and even equi-uniformly continuous $\!)$ on $K$.
\end{lemma}
{\sc Proof}:
% Note first that each single $\widehat{\mu}_n\in C(H,\CC)$ is uniformly 
% continuous because $\mu_n\in\cP(W)\subset\cP(H)$, and therefore 
% $\|\widehat{\mu}_n\|_{\infty}\le\|\mu_n\|=1$.
Let $\varepsilon>0$ and $V_{\varepsilon} :=\{ k\in\widehat{H}\mid 
\sup_{x\in W} |\langle k,x\rangle -1| < \varepsilon \}$. Note that
$0\in V_{\varepsilon}$ (since the trivial character is 
$\langle 0,x\rangle\equiv 1$) and that $V_{\varepsilon}$ is a neighbourhood 
of $0$ in $\widehat{H}$ because it is a typical 
open set in the compact-open topology of $\widehat{H}\subset C(H,\CC)$.
Fix some $k_1\in K$ and choose $k_2\in K$ so that $k_1-k_2\in V_{\varepsilon}$.
Then we have, for all $x\in W$,
$$\Bigl|\,\overline{\langle k_1,x\rangle} - \overline{\langle k_2,x\rangle}
\,\Bigr|\; = \; \bigl|\langle k_1-k_2,x\rangle -1 \bigr| \; < \varepsilon\, , $$
where we used $|\langle k_2,x\rangle|=1$ and the multiplication rule
for characters.

Let $\widehat{\mu}_n$ be an arbitrary element of our sequence. We then get
\begin{eqnarray*}
   \bigl|\widehat{\mu}_n(k_1) - \widehat{\mu}_n(k_2)\bigr| & = &
   \left|\int_H \Bigl(\overline{\langle k_1,x\rangle} - 
   \overline{\langle k_2,x\rangle}\Bigr) \, {\rm d}\mu_n(x) \right| \\
   & \le & \int_W \Bigl|\,\overline{\langle k_1,x\rangle} - 
   \overline{\langle k_2,x\rangle}\,\Bigr| \, {\rm d}\mu_n(x) \\
   & < & \varepsilon \, \|\mu_n\| \; = \; \varepsilon \, .
\end{eqnarray*}
Since this is independent both of $n$ and of $k_1\in K$, equi-uniform
continuity of $\{\widehat{\mu}_n\}$ on $K$ follows. \qed

\smallskip

{\sc Remark}: With little extra complication, the result can 
be extended to vaguely convergent sequences of measures
$\mu_n\in\cM_+^m(H)$, see \cite[Lemma 23.7]{Bauer} for a proof that can
easily be adapted to this case.

\begin{theorem} {\rm (Continuity Theorem of P.\ L\'evy)}\label{levy}
   $\,$Let $W,H,\widehat{H}$ be defined as in Lemma\ {\rm \ref{equi-cont}}.
   If $(\mu_n)_{n\in\NN}$ is a sequence of measures in $\cP(W)$ that
   vaguely converges to $\mu\in\cP(W)$, then the corresponding sequence
   of Fourier-Stieltjes transforms, $(\widehat{\mu}_n)_{n\in\NN}$,
   converges compactly to $\widehat{\mu}$.
\end{theorem} \index{L\'evy's continuity theorem}
{\sc Proof}:
Each of the functions $x\mapsto\langle k,x\rangle$, $k\in\widehat{H}$,
lies in $C(H,\CC)$. Therefore, 
$\lim_{n\to\infty} \widehat{\mu}_n(k) = \widehat{\mu}(k)$ by the very
definition of vague convergence. So, we have a sequence of pointwise
converging functions in $C(H,\CC)$ that are equi-continuous
on compact subsets $K\subset\widehat{H}$ due to Lemma~\ref{equi-cont}, 
hence the convergence of $(\widehat{\mu}_n)_{n\in\NN}$ is uniform on $K$
by Lemma~\ref{point-compact} which proves the assertion. \qed

\smallskip

{\sc Remark}: Once again, the restriction to $\cP(W)$ is not essential,
and the result holds for sequences in $\cM_+^m(H)$ as well, compare
\cite[Thm.\ 23.8]{Bauer}, but we do not need the stronger result here.

{}Finally, we formulate the approximation property that we need in
Section \ref{model-weyl} to prove Theorem~\ref{weyl-thm}.
\begin{lemma}\label{step}
  Let $H$ be an LCAG with Haar measure $\mu$, and let
  $W$ be a Riemann measurable compact subset of $H$.
  Then each $f\in C(H,\CC)$ with ${\rm supp}(f)\subset W$
  can be uniformly approximated by
  a sequence of step functions $\psi^{(\ell)}$ of the form
$$ \psi^{(\ell)} \; = \; \sum_{i=1}^{n(\ell)} 
   c^{(\ell)}_{i} \, \bs{1}^{}_{U_i} $$
  where $c^{(\ell)}_{i}\in\CC$ and $W=\bigcup_{i=1}^{n(\ell)} U_i$ is
  a partition into pairwise disjoint sets $U_i\subset W$ that are
  all Riemann measurable.
\end{lemma}
{\sc Proof}:
Although this is a standard type of result, we give an explicit proof
because the additional condition of Riemann measurability of the partition
requires some attention. 
It is sufficient to prove the statement for real functions
because any $f$ may be split into its real and imaginary parts, both of 
which are continuous and supported on $W$.

Assume $f$ is real.
{}For all  $s\in \RR $, $f^{-1}(s)$ is a closed subset of $H$, hence 
measurable.  {}For $0\neq s\neq t\neq 0$, 
$f^{-1}(s)$ and $f^{-1}(t)$ are disjoint subsets of $W$, and we thus have
$\sum_{s\in\RR\setminus\{0\}}\mu(f^{-1}(s))\le \mu(W) < \infty$.
But this means, as usual, that 
$$ P \; := \; \{ s\in\RR\mid \mu(f^{-1}(s)) > 0\} $$
is at most a {\em countable} set, called the set of ``bad'' points.
Since $f$ is continuous, $f(W)\subset\RR$ is compact, hence
$f(W)\subset\ (a,b)$ for some $a,b\in\RR$.

{}Fix some $\varepsilon>0$ and choose an integer $n>(b-a)/\varepsilon$.
We can then subdivide $[a,b)$ into non-empty intervals $I^{}_1,\dots,I^{}_n$, 
each of the form $[\alpha_i,\beta_i)$ with $\alpha_i,\beta_i\not\in P$ and
$0<{\rm length}(I_i)=\beta_i-\alpha_i<\varepsilon$. Define 
$U_i=f^{-1}(I_i)\cap W$.
Then $\{U_1,\dots,U_n\}$ clearly is a partition of $W$. In addition,
we have, for all $1\le i\le n$,
$$ \left. \begin{array}{r}
          f^{-1}((\alpha_i,\beta_i)) \ \mbox{is open} 
          \vphantom{\prod^{1}_{2}}\\
          f^{-1}([\alpha_i,\beta_i]) \ \mbox{is closed}
          \vphantom{\prod^{1}_{2}}\\
          \mu\{f^{-1}(\alpha_i)\cup f^{-1}(\beta_i)\}=0
          \vphantom{\prod^{1}_{2}}
   \end{array} \right\} 
   \; \Longrightarrow \; \mbox{each } U_i \mbox{ is Riemann measurable}.
$$

Now, $\overline{U}_i\subset W$ is closed and hence compact. Define
$c_i:=\inf\{f(x)\mid x\in U_i\}$ and 
$\psi:=\sum_{i=1}^{n} c_i \,\bs{1}^{}_{U_i}$. Then $\psi$ is supported on $W$.
If $x\in W$, then $x$ is in precisely one of the sets of the partition,
$U_i$ say, and $|f(x) - \psi(x)| = |f(x) - c_i| \le {\rm length}(I_i)
<\varepsilon$, independently of $x$. Consequently,
$\|f-\psi\|_{\infty}<\varepsilon$. {}Finally, we can construct
a sequence $\psi^{(\ell)}$ this way (e.g.\ via $\varepsilon^{(\ell)}=1/\ell$)
which establishes our claim. \qed

\smallskip

{\sc Remark}: The result of this Lemma extends to all continuous functions
on $W$, even if they do not define continuous functions on $H$. This
follows from Tietze's extension theorem \cite[Thm.\ IV.11]{RS}.

\bibliographystyle{amsalpha}

\begin{thebibliography}{WWA}

\bibitem[A]{Akin} 
E.~Akin, 
\textit{The General Topology of Dynamical Systems}, 
% Graduate Studies in Mathematics, Vol. 1,
Amer.\ Math.\ Soc., Providence, RI (1993);
corr.\ reprint (1996).

\bibitem[AGS]{AGS}
R.~Ammann, B.~Gr\"unbaum and G.\thinspace C.~Shephard,
\textit{Aperiodic Tiles}, 
Discr.\ Comput.\ Geom.\ \textbf{8} (1992), 1--25.

\bibitem[B]{MB}
M.~Baake,
\textit{A guide to mathematical quasicrystals},
in: \textit{Quasicrystals}, eds.\ J.-B.\ Suck,
M.\ Schreiber and P.\ H\"au{\ss}ler, Springer,
Berlin (2001), in press; preprint math-ph/9901014.

\bibitem[BJ]{BJ}
M.~Baake and D.~Joseph,
\textit{Ideal and defective vertex configurations in the planar
octagonal quasilattice}, Phys.\ Rev.\ \textbf{B 42} (1990), 8091--8102.

\bibitem[BM1]{BM}
M.~Baake and R.\thinspace V.~Moody, 
\textit{Self-similarities and invariant
densities for model sets}, in: \textit{Algebraic Methods and
Theoretical Physics}, ed.\ Y.\ Saint-Aubin, 
Springer, New York (2000), in press;
preprint math-ph/9809006.

\bibitem[BM2]{BMmulti}
M.~Baake and R.\thinspace V.~Moody, 
\textit{Multi-component model sets and invariant densities}, 
in: \textit{Aperiodic '97},
eds.\ M.~de Boissieu, J.-L.~Verger-Gaugry and R.~Currat, 
World Scientific, Singapore (1998), pp.\ 9--20;
math-ph/9809005.

\bibitem[BMS]{BMS}
M.~Baake, R.\thinspace V.~Moody and M.~Schlottmann,
\textit{Limit-(quasi)periodic point sets as model sets with
$p$-adic internal spaces},
J.\ Phys.\ \textbf{A 31} (1998), 5755--5765;
math-ph/9901008.

\bibitem[Bau]{Bauer}
H.~Bauer, 
\textit{Probability Theory}, 
% Studies in Mathematics, Vol.\ 23,
de Gruyter, Berlin (1996).

\bibitem[BoGi]{BG}
J.~M.~Borwein and R.~Girgensohn,
\textit{Functional equations and distribution functions},
Results in Math.\ \textbf{26} (1994), 229--237.

\bibitem[Bou]{Bourbaki}
N.~Bourbaki, 
\textit{General Topology}, vol.~1, Addison-Wesley, Reading, MA (1966).

\bibitem[D]{Dieu}
J.~Dieudonn\'e,
\textit{Treatise on Analysis}, vol.\ II, Academic Press, New York (1970).

\bibitem[HRB]{HRB}
J.~Hermisson, C.~Richard and M.~Baake,
\textit{A guide to the symmetry structure of quasiperiodic tiling
classes}, J.\ Phys.\ I (France) \textbf{7} (1997), 1003--1018.

\bibitem[HeRo]{HR}
E.~Hewitt and K.\thinspace A.~Ross,
\textit{Abstract Harmonic Analysis}, vol.\ 1,
2nd ed., Springer, New York (1979);
corr.\ 3rd printing (1997).

\bibitem[HM]{HM}
M.~H\"offe and R.\thinspace V.~Moody, 
\textit{Self-similar measures for $p$-adic model sets},
in preparation.

\bibitem[Hof1]{Hof-diff}
A.~Hof, 
\textit{On diffraction by aperiodic structures},
Commun.\ Math.\ Phys.\ \textbf{169} (1995), 25--43.

\bibitem[Hof2]{Hof}
A.~Hof, 
\textit{Uniform distribution and the projection method}, 
in: \textit{Quasicrystals and Discrete Geometry},
ed.\ J.~Patera, Fields Institute Monographs, vol.\ 10, 
Amer.\ Math.\ Soc., Providence, RI (1998), pp.\ 201--206.

\bibitem[Hut]{Hutchinson}
J.\thinspace E.~Hutchinson,  
\textit{Fractals and self-similarity}, 
Indiana Univ.\ Math.\ J.\ \textbf{30} (1981), 713--743.

% \bibitem[JW]{JW}
% B.~Jessen and A.~Wintner,  
% \textit{Distribution functions and the Riemann zeta function}, 
% Trans.\ Amer.\ Math.\ Soc.\ \textbf{38} (1935), 48--88.

\bibitem[JLS]{JLS}
R.\thinspace Q.~Jia, S.\thinspace L.~Lee and A.~Sharma, 
\textit{Spectral properties of continuous refinement operators},
Proc.\ Amer.\ Math.\ Soc.\ \textbf{126} (1998), 729--737.

\bibitem[JL]{JL}
Q.~Jiang and S.\thinspace L.~Lee,
\textit{Spectral properties of matrix continuous refinement operators},
Adv.\ Comput.\ Math.\ \textbf{7} (1997), 361--382.

\bibitem[KT]{KT}
S.~Karlin and H.\thinspace M.~Taylor,
\textit{A First Course in Stochastic Processes},
2nd ed., Academic Press, Boston, MA (1975)

\bibitem[Kel]{Kelley} 
J.\thinspace L.~Kelley, 
\textit{General Topology}, 
Van Nostrand, Princeton, NJ (1955);
reprint, Springer, New York (1975).

\bibitem[KN]{KN}
L.~Kuipers and H.~Niederreiter,
\textit{Uniform Distribution of Sequences},
Wiley, New York (1974).

\bibitem[LW]{LW}
J.\thinspace C.~Lagarias and Y.~Wang, 
\textit{Substitution Delone sets},
preprint (2000). 

\bibitem[Lal]{Lal}
S.~P.~Lalley,
\textit{Random Series in powers of algebraic integers:
Hausdorff dimension of the limit distribution},
J.\ London Math.\ Soc.\ (2) \textbf{57} (1998), 629--654.

\bibitem[LM]{LM}
J.-Y.~Lee and R.\thinspace V.~Moody, 
\textit{Lattice substitution systems and model sets},
Discr.\ Comput.\ Geom., to appear; 
preprint math.MG/0002019.

\bibitem[Mey]{Meyer}
Y.~Meyer,
\textit{Algebraic Numbers and Harmonic Analysis},
North-Holland, Amsterdam (1972).

\bibitem[M]{Moody}
R.\thinspace V.~Moody,
\textit{Meyer sets and their duals}, in: 
\textit{The Mathematics of Long-Range Aperiodic Order},
ed.\ R.\thinspace V.\ Moody, NATO ASI Series C 489,
Kluwer, Dordrecht (1997), pp.\ 403--441.

\bibitem[P]{P}
P.\thinspace A.\thinspace B.~Pleasants,
\textit{Designer quasicrystals: cut-and-project sets with
pre-assigned properties}, 
in: \textit{Directions in Mathematical Quasicrystals},
eds.\ M.~Baake and R.\thinspace V.~Moody, CRM Monograph Series, 
Amer.\ Math.\ Soc., Providence, RI (2000); this volume.

\bibitem[Q]{Q}
B.~von Querenburg,
\textit{Mengentheoretische Topologie},
2nd ed., Springer, Berlin (1979).

\bibitem[RS]{RS}
M.~Reed and B.~Simon,
\textit{Methods of Modern Mathematical Physics I: Functional Analysis},
2nd ed., Academic Press, San Diego, CA (1980).

\bibitem[Ru]{Rudin}
W.~Rudin,
\textit{Fourier Analysis on Groups},
Wiley, New York (1962); reprint (1990).

\bibitem[Sch1]{Martin}
M.~Schlottmann, 
\textit{Cut-and-project sets in locally compact Abelian groups}, 
in: \textit{Quasicrystals and Discrete Geometry}, 
ed.\ J.~Patera, Fields Institute Monographs, vol.\ 10, 
Amer.\ Math.\ Soc., Providence, RI (1998), 247--264.

\bibitem[Sch2]{Martin2}
M.~Schlottmann, 
\textit{Generalized model sets and dynamical systems},
in: \textit{Directions in Mathematical Quasicrystals},
eds.\ M.~Baake and R.\thinspace V.~Moody, CRM Monograph series, 
Amer.\ Math.\ Soc., Providence, RI (2000); this volume.

\bibitem[So]{Sol}
B.~ Solomyak, 
\textit{On the random series $\sum \pm \lambda^n$ $($an Erd\H{o}s problem$\,)$}, 
Ann.\ Math.\ \textbf{142} (1995), 611--625.

\bibitem[We]{Weyl}
H.~Weyl, 
\textit{\"Uber die Gleichungverteilung von Zahlen mod.\ Eins}, 
Math.\ Ann.\ \textbf{77} (1916), 313--352.

\bibitem[Wi]{Wicks}
K.\thinspace R.~Wicks, 
\textit{Fractals and Hyperspaces}, 
Lecture Notes in Mathematics, vol.\ 1492, 
Springer, Berlin (1991).


\end{thebibliography}

\printindex
\end{document}